\documentclass[12pt,reqno]{article}
\usepackage[usenames]{color}
\usepackage{amssymb}
\usepackage{graphicx}
\usepackage{amscd}

\usepackage[colorlinks=true,
linkcolor=webgreen,
filecolor=webbrown,
citecolor=webgreen]{hyperref}

\definecolor{webgreen}{rgb}{0,.5,0}
\definecolor{webbrown}{rgb}{.6,0,0}

\usepackage{color}
\usepackage{fullpage}
\usepackage{float}

\usepackage{graphics,amsmath,amssymb}
\usepackage{amsthm}
\usepackage{amsfonts}
\usepackage{latexsym}
\usepackage{epsf}

\usepackage{verbatim}
\usepackage{enumerate}

\usepackage{tikz}
\usepackage{forest}

\usepackage{mathtools}

\usepackage{array}

\usepackage[blocks]{authblk}

\theoremstyle{plain}
\newtheorem{theorem}{Theorem}
\newtheorem{corollary}[theorem]{Corollary}

\newtheorem{proposition}[theorem]{Proposition}

\theoremstyle{definition}

\theoremstyle{remark}

\newcommand{\PrefSquareFourSmallEditable}[0]{

\draw[black, ultra thick] (-3,3) -- (3,3);
\draw[black, ultra thick] (-3,3) -- (-3,-3);
\draw[black, ultra thick] (3,3) -- (3,-3);
\draw[black, ultra thick] (-3,-3) -- (3,-3);
\draw[black, ultra thick] (0,3) -- (0,-3);
\draw[black, ultra thick] (3,0) -- (-3,0);
\draw[black, ultra thin] (-2.25,3) -- (-2.25,-3);
\draw[black, ultra thick] (-1.5,3) -- (-1.5,-3);
\draw[black, ultra thin] (-0.75,3) -- (-0.75,-3);
\draw[black, ultra thin] (0.75,3) -- (0.75,-3);
\draw[black, ultra thick] (1.5,3) -- (1.5,-3);
\draw[black, ultra thin] (2.25,3) -- (2.25,-3);
\draw[black, ultra thin] (3,-2.25) -- (-3,-2.25);
\draw[black, ultra thick] (3,-1.5) -- (-3,-1.5);
\draw[black, ultra thin] (3,-0.75) -- (-3,-0.75);
\draw[black, ultra thin] (3,0.75) -- (-3,0.75);
\draw[black, ultra thick] (3, 1.5) -- (-3,1.5);
\draw[black, ultra thin] (3,2.25) -- (-3,2.25);
\draw[black, ultra thin] (-2.625,3) -- (-2.625,-3);
\draw[black, ultra thin] (-1.875,3) -- (-1.875,-3);
\draw[black, ultra thin] (-1.125,3) -- (-1.125,-3);
\draw[black, ultra thin] (-0.375,3) -- (-0.375,-3);
\draw[black, ultra thin] (0.375,3) -- (0.375,-3);
\draw[black, ultra thin] (1.125,3) -- (1.125,-3);
\draw[black, ultra thin] (1.875,3) -- (1.875,-3);
\draw[black, ultra thin] (2.625,3) -- (2.625,-3);
\draw[black, ultra thin] (-3,-2.625) -- (3,-2.625);
\draw[black, ultra thin] (-3,-1.875) -- (3,-1.875);
\draw[black, ultra thin] (-3,-1.125) -- (3,-1.125);
\draw[black, ultra thin] (-3,-0.375) -- (3,-0.375);
\draw[black, ultra thin] (-3,2.625) -- (3,2.625);
\draw[black, ultra thin] (-3,1.875) -- (3,1.875);
\draw[black, ultra thin] (-3,1.125) -- (3,1.125);
\draw[black, ultra thin] (-3,0.375) -- (3,0.375);

\node at (-2.45,2.8) [text width=1cm] {\tiny 1};
\node at (-2.45,2.425) [text width=1cm] {\tiny 2};
\node at (-2.45,2.05) [text width=1cm] {\tiny 3};
\node at (-2.45,1.675) [text width=1cm] {\tiny 4};
\node at (-2.45,1.3) [text width=1cm] {\tiny 1};
\node at (-2.45,0.925) [text width=1cm] {\tiny 2};
\node at (-2.45,0.55) [text width=1cm] {\tiny 3};
\node at (-2.45,0.175) [text width=1cm] {\tiny 4};
\node at (-2.45,-0.2) [text width=1cm] {\tiny 1};
\node at (-2.45,-0.575) [text width=1cm] {\tiny 2};
\node at (-2.45,-0.95) [text width=1cm] {\tiny 3};
\node at (-2.45,-1.325) [text width=1cm] {\tiny 4};
\node at (-2.45,-1.7) [text width=1cm] {\tiny 1};
\node at (-2.45,-2.075) [text width=1cm] {\tiny 2};
\node at (-2.45,-2.45) [text width=1cm] {\tiny 3};
\node at (-2.45,-2.825) [text width=1cm] {\tiny 4};
\node at (-2.025,3.25) [text width=1cm] {\tiny 1};
\node at (-1.65,3.25) [text width=1cm] {\tiny 2};
\node at (-1.275,3.25) [text width=1cm] {\tiny 3};
\node at (-0.9,3.25) [text width=1cm] {\tiny 4};
\node at (-0.525,3.25) [text width=1cm] {\tiny 1};
\node at (-0.15,3.25) [text width=1cm] {\tiny 2};
\node at (0.225,3.25) [text width=1cm] {\tiny 3};
\node at (0.6,3.25) [text width=1cm] {\tiny 4};
\node at (0.99,3.25) [text width=1cm] {\tiny 1};
\node at (1.35,3.25) [text width=1cm] {\tiny 2};
\node at (1.725,3.25) [text width=1cm] {\tiny 3};
\node at (2.1,3.25) [text width=1cm] {\tiny 4};
\node at (2.485,3.25) [text width=1cm] {\tiny 1};
\node at (2.855,3.25) [text width=1cm] {\tiny 2};
\node at (3.225,3.25) [text width=1cm] {\tiny 3};
\node at (3.6,3.25) [text width=1cm] {\tiny 4};

\node at (-1.575, 3.8) [text width=1cm] {\normalsize A};
\node at (0.025, 3.8) [text width=1cm] {\normalsize B};
\node at (1.325, 3.8) [text width=1cm] {\normalsize C};
\node at (2.825, 3.8) [text width=1cm] {\normalsize D};
\node at (-3.35, 2.25) [text width=1cm] {\normalsize E};
\node at (-3.35, 0.725) [text width=1cm] {\normalsize F};
\node at (-3.35, -0.725) [text width=1cm] {\normalsize G};
\node at (-3.35, -2.25) [text width=1cm] {\normalsize H};
}

\newcommand{\DoubleEdit}[0]{

\draw[black, ultra thick] (2,2) -- (-2,2);
\draw[black, ultra thick] (-2,2) -- (-2,-2);
\draw[black, ultra thick] (-2,-2) -- (2,-2);
\draw[black, ultra thick] (2,2) -- (2,-2);
\draw[black, very thick] (0,2) -- (0,-2);
\draw[black, very thick] (-2,0) -- (2,0);
\draw[black, very thin] (-1,2) -- (-1,-2);
\draw[black, very thin] (1,2) -- (1,-2);
\draw[black, very thin] (2,1) -- (-2,1);
\draw[black, very thin] (2,-1) -- (-2,-1);
}

\newcommand{\TripleEdit}[0]{
\draw[black, ultra thick] (3,3) -- (-3,3);
\draw[black, ultra thick] (-3,3) -- (-3,-3);
\draw[black, ultra thick] (-3,-3) -- (3,-3);
\draw[black, ultra thick] (3,-3) -- (3,3);
\draw[black, very thick] (1,3) -- (1,-3);
\draw[black, very thick] (-1,3) -- (-1,-3);
\draw[black, very thick] (-3,1) -- (3,1);
\draw[black, very thick] (-3,-1) -- (3,-1);
\draw[black, very thin] (-7/3, 3) -- (-7/3,-3);
\draw[black, very thin] (-5/3, 3) -- (-5/3,-3);
\draw[black, very thin] (-1/3, 3) -- (-1/3,-3);
\draw[black, very thin] (1/3, 3) -- (1/3,-3);
\draw[black, very thin] (5/3, 3) -- (5/3,-3);
\draw[black, very thin] (7/3, 3) -- (7/3,-3);
\draw[black, very thin] (3, -7/3) -- (-3,-7/3);
\draw[black, very thin] (3, -5/3) -- (-3,-5/3);
\draw[black, very thin] (3, -1/3) -- (-3,-1/3);
\draw[black, very thin] (3, 1/3) -- (-3,1/3);
\draw[black, very thin] (3, 5/3) -- (-3,5/3);
\draw[black, very thin] (3, 7/3) -- (-3,7/3);
}
\newcommand{\FirstColumnTripleEdit}[9]{
\node at (-8/3,8/3) {\normalsize #1 };
\node at (-8/3,2) {\normalsize #2 };
\node at (-8/3,4/3) {\normalsize #3 };
\node at (-8/3,2/3) {\normalsize #4 };
\node at (-8/3,0) {\normalsize #5 };
\node at (-8/3,-2/3) {\normalsize #6 };
\node at (-8/3,-4/3) {\normalsize #7 };
\node at (-8/3,-6/3) {\normalsize #8 };
\node at (-8/3,-8/3) {\normalsize #9 };
}

\newcommand{\SecondColumnTripleEdit}[9]{
\node at (-6/3,8/3) {\normalsize #1 };
\node at (-6/3,2) {\normalsize #2 };
\node at (-6/3,4/3) {\normalsize #3 };
\node at (-6/3,2/3) {\normalsize #4 };
\node at (-6/3,0) {\normalsize #5 };
\node at (-6/3,-2/3) {\normalsize #6 };
\node at (-6/3,-4/3) {\normalsize #7 };
\node at (-6/3,-6/3) {\normalsize #8 };
\node at (-6/3,-8/3) {\normalsize #9 };
}
\newcommand{\ThirdColumnTripleEdit}[9]{
\node at (-4/3,8/3) {\normalsize #1 };
\node at (-4/3,2) {\normalsize #2 };
\node at (-4/3,4/3) {\normalsize #3 };
\node at (-4/3,2/3) {\normalsize #4 };
\node at (-4/3,0) {\normalsize #5 };
\node at (-4/3,-2/3) {\normalsize #6 };
\node at (-4/3,-4/3) {\normalsize #7 };
\node at (-4/3,-6/3) {\normalsize #8 };
\node at (-4/3,-8/3) {\normalsize #9 };
}
\newcommand{\FourthColumnTripleEdit}[9]{
\node at (-2/3,8/3) {\normalsize #1 };
\node at (-2/3,2) {\normalsize #2 };
\node at (-2/3,4/3) {\normalsize #3 };
\node at (-2/3,2/3) {\normalsize #4 };
\node at (-2/3,0) {\normalsize #5 };
\node at (-2/3,-2/3) {\normalsize #6 };
\node at (-2/3,-4/3) {\normalsize #7 };
\node at (-2/3,-6/3) {\normalsize #8 };
\node at (-2/3,-8/3) {\normalsize #9 };
}
\newcommand{\FifthColumnTripleEdit}[9]{
\node at (0/3,8/3) {\normalsize #1 };
\node at (0/3,2) {\normalsize #2 };
\node at (0/3,4/3) {\normalsize #3 };
\node at (0/3,2/3) {\normalsize #4 };
\node at (0/3,0) {\normalsize #5 };
\node at (0/3,-2/3) {\normalsize #6 };
\node at (0/3,-4/3) {\normalsize #7 };
\node at (0/3,-6/3) {\normalsize #8 };
\node at (0/3,-8/3) {\normalsize #9 };
}
\newcommand{\SixthColumnTripleEdit}[9]{
\node at (2/3,8/3) {\normalsize #1 };
\node at (2/3,2) {\normalsize #2 };
\node at (2/3,4/3) {\normalsize #3 };
\node at (2/3,2/3) {\normalsize #4 };
\node at (2/3,0) {\normalsize #5 };
\node at (2/3,-2/3) {\normalsize #6 };
\node at (2/3,-4/3) {\normalsize #7 };
\node at (2/3,-6/3) {\normalsize #8 };
\node at (2/3,-8/3) {\normalsize #9 };
}
\newcommand{\SeventhColumnTripleEdit}[9]{
\node at (4/3,8/3) {\normalsize #1 };
\node at (4/3,2) {\normalsize #2 };
\node at (4/3,4/3) {\normalsize #3 };
\node at (4/3,2/3) {\normalsize #4 };
\node at (4/3,0) {\normalsize #5 };
\node at (4/3,-2/3) {\normalsize #6 };
\node at (4/3,-4/3) {\normalsize #7 };
\node at (4/3,-6/3) {\normalsize #8 };
\node at (4/3,-8/3) {\normalsize #9 };
}
\newcommand{\EighthColumnTripleEdit}[9]{
\node at (6/3,8/3) {\normalsize #1 };
\node at (6/3,2) {\normalsize #2 };
\node at (6/3,4/3) {\normalsize #3 };
\node at (6/3,2/3) {\normalsize #4 };
\node at (6/3,0) {\normalsize #5 };
\node at (6/3,-2/3) {\normalsize #6 };
\node at (6/3,-4/3) {\normalsize #7 };
\node at (6/3,-6/3) {\normalsize #8 };
\node at (6/3,-8/3) {\normalsize #9 };
}
\newcommand{\NinthColumnTripleEdit}[9]{
\node at (8/3,8/3) {\normalsize #1 };
\node at (8/3,2) {\normalsize #2 };
\node at (8/3,4/3) {\normalsize #3 };
\node at (8/3,2/3) {\normalsize #4 };
\node at (8/3,0) {\normalsize #5 };
\node at (8/3,-2/3) {\normalsize #6 };
\node at (8/3,-4/3) {\normalsize #7 };
\node at (8/3,-6/3) {\normalsize #8 };
\node at (8/3,-8/3) {\normalsize #9 };
}

\newcommand{\TripleSquareEdit}[0]{
\draw[black, ultra thick] (3,3) -- (-3,3);
\draw[black, ultra thick] (-3,3) -- (-3,-3);
\draw[black, ultra thick] (-3,-3) -- (3,-3);
\draw[black, ultra thick] (3,-3) -- (3,3);
\draw[black, thick] (1,3) -- (1,-3);
\draw[black, thick] (-1,3) -- (-1,-3);
\draw[black, thick] (-3,1) -- (3,1);
\draw[black, thick] (-3,-1) -- (3,-1);
}

\begin{document}

\title{The Stable Matching Problem and Sudoku}
\author{Matvey Borodin}
\author{Eric Chen}
\author{Aidan Duncan}
\author{Boyan Litchev}
\author{Jiahe Liu}
\author{Veronika Moroz}
\author{Matthew Qian}
\author{Rohith Raghavan}
\author{Garima Rastogi}
\author{Michael Voigt}
\affil{PRIMES STEP}
\author{Tanya Khovanova}
\affil{MIT}

\maketitle

\begin{abstract}
Are you having trouble getting married? These days, there are lots of products on the market for dating, from apps to websites and matchmakers, but we know a simpler way! That's right --- your path to coupled life isn't through Tinder: it's through Sudoku! Read our fabulous paper where we explore the Stable Marriage Problem to help you find happiness and stability in marriage through math. As a bonus, you get two Sudoku puzzles with a new flavor.
\end{abstract}

\section{Introduction}\label{sec:intro}

The stable marriage problem (SMP) is the problem of matching $n$ men and $n$ women into married couples depending on their preferences so that the matching is stable. To make it easier to define a stable matching, we first define an unstable matching. An unstable matching is a matching such that there exist two people who prefer each other to their partners. Thus, a stable matching is a matching where no such pairs of people exist.

In a famous paper dating back to 1962, Gale and Shapley \cite{GS1962} suggested an algorithm, now called the Gale-Shapley algorithm, that always finds a stable matching. Theoretically, it is possible to have distinct stable matchings for the same set of preferences. The original Gale-Shapley algorithm is based on men proposing to women. It favors men; that is, the result is man-optimal and woman-pessimal.

The stable marriage problem is extremely versatile as it can be applied to matching applicants and jobs, users and servers, and so on. Because of its numerous applications, it attracts a lot of researchers. But now, it is time to turn to Sudoku.

A Sudoku is a logic puzzle played on a partially-filled 9 by 9 grid. The objective is to complete the grid by filling it with the integers 1 through 9 so that every row, column, and 3 by 3 block contains distinct digits. Note that the blocks are the nine non-overlapping 3 by 3 subgrids. By tradition, a Sudoku puzzle has a unique solution. The mathematics behind Sudoku is covered in an awesome book by Rosenhouse and Taalman \cite{RT2011}.

How is the stable marriage problem related to Sudoku? This is what this paper is about. 

One of the main goals of the paper is to draw parallels between SMP and Sudoku. While doing this, we invented a new type of Sudoku, which we call a joint-groups Sudoku. We describe this Sudoku type in detail and also discuss the corresponding preference profiles for men and women.

Now we provide a more detailed description of the paper.

Specifically, Section~\ref{sec:defs} is devoted to definitions and background related to the stable marriage problem and Sudokus.

Section~\ref{sec:SMP2Sudokus} translates SMP into the language of Sudoku. We explain how to draw a preference profile in a Sudoku grid. We also show how to find blocking pairs and explain the Gale-Shapley algorithm in terms of Sudoku.

Section~\ref{sec:Sudokus2SMP} translates Sudokus into the language of SMP. We define non-overlapping preference profiles, which are profiles that could correspond to different digits in a complete Sudoku grid. We explain what it means to solve a Sudoku puzzle for these preference profiles. 

We start Section~\ref{sec:specialprofiles} with defining pseudo-Latin profiles introduced in \cite{T2002}. We generalize this notion and introduce mutually Latin profiles. These are the profiles where the men's and the women's preference matrices are Latin squares. We also introduce disjoint profiles inspired by disjoint-groups Sudokus. 

Section~\ref{sec:joint} introduces a new Sudoku type called joint-groups Sudoku. The profiles corresponding to this Sudoku are called joint profiles. We show that a joint profile is uniquely defined by a key function. We also show that joint profiles are mutually Latin profiles. In addition, we classify joint profiles for $n = 3$, and discuss some properties of their stable matchings. Finally, we show that joint-groups Sudokus exist and give two joint-groups Sudoku puzzles for readers to solve.

In Section~\ref{sec:n=2} we apply all the theory we developed to completely analyze the case of $n=2$.

Section~\ref{sec:answers} contains answers to the puzzles from Section~\ref{sec:joint}.

\section{Definitions}\label{sec:defs}

In the stable matching problem we study in this paper, we have $n$ men and $n$ women who prefer getting married over being single.  As is traditional in this problem, we assume that men only marry women and vice versa. Each person ranks people of the other gender without ties. We call these $2n$ sets of preferences a \textit{preference profile}.

Once the preference profile is set, people are matched into marriages. We call the set of all these marriages \textit{unstable} if there is a man $M$ and woman $W$ such that they prefer each other to their spouses, and \textit{stable} if no such man/woman pair exists. In an unstable set of marriages, the pair $M$ and $W$ described above is called a \textit{blocking pair} or a \textit{rogue couple}. By definition, a set of marriages without a blocking pair is stable.

There exist algorithms to find a stable matching given a preference profile. The first and most famous is the Gale-Shapley algorithm \cite{GS1962} (also known as the deferred acceptance algorithm), which involves a number of rounds:
\begin{itemize}
\item In the first round, each man proposes to the woman he prefers most, and then each woman replies ``maybe'' to her most preferred suitor and ``no'' to all other suitors. She is then provisionally ``engaged'' to the suitor she most prefers so far, and that suitor is likewise provisionally engaged to her.
\item In each subsequent round, first each unengaged man proposes to the most-preferred woman to whom he has not yet proposed (regardless of whether the woman is already engaged), and then each woman replies ``maybe'' if she is currently not engaged or if she prefers this man over her current provisional partner (in this case, she rejects her current provisional partner who becomes unengaged).
\item This process is repeated until everyone is engaged. At this point, all engaged couples form a stable matching. They can start planning their weddings.
\end{itemize}

For $n$ men and $n$ women, the matching algorithm takes at most $n^2-n+1$ rounds and terminates in a stable matching \cite{GS1962}.

Let $S$ be the set of all possible stable matchings. We call man $M$ \textit{a valid partner} of woman $W$ if there exists some stable matching $s \in S$ where they are matched. A matching is called \textit{man-optimal} if each man receives his best valid partner, and \textit{man-pessimal} if each man receives his worst valid partner.

The Gale-Shapley algorithm, where men are proposing, results in a man-optimal stable matching \cite{GI1989}. If a matching is man-optimal, then it is also woman-pessimal \cite{GI1989}. It follows that if the men-proposing and women-proposing algorithms end in the same matching, then this is the only stable matching.

\subsection{Latin Stable Marriage Problem}

The \textit{Latin} marriage problem is a subset of the stable marriage
problem such that the sum of the mutual rankings for a man and a woman is $n+1$ \cite{T2002}. This means the matrix of the men's preferences and the matrix of the women's preferences each form a Latin square \cite{T2002}. This setup is sometimes called a pseudo-Latin stable matching. We call the profiles in a pseudo-Latin stable matching \textit{pseudo-Latin} profiles.  We have to point out that there exist preference profiles that form Latin squares for men and women but are not pseudo-Latin profiles. 

The Latin marriage problem is interesting because pseudo-Latin profiles tend to produce a lot of stable matchings. Sequence A069124 in the OEIS \cite{OEIS} describes the number of possible stable matchings in a pseudo-Latin profile with $n$ men and $n$ women. It starts as
\[1,\ 2,\ 3,\ 10,\ 12,\ 32,\ 42,\ 268,\ 288,\ 656,\ 924,\ 4360,\ 3816,\ 11336,\ \ldots.\]

This sequence counts the number of stable matchings for special types of profiles, thus providing a lower bound in the number of stable matchings for any profile.

\subsection{Preference profiles}

The total number of different preference profiles for $n$ men and $n$ women is 
\[(n!)^{2n}\]
since each person of one gender can rank all the people of the opposite gender in $n!$ ways. The corresponding sequence is sequence A185141 in the OEIS \cite{OEIS}. The sequence starts as follows, where the first term corresponds to $n=1$:
\[1,\ 16,\ 46656,\  110075314176,\ 619173642240000000000, \ldots.\]

%	1, 1, 16, 46656, 110075314176, 619173642240000000000, 19408409961765342806016000000000000, 6823819180249038753817675898369448345600000000000000, 48789725533845219197010193096946682961355723304326670581760000000000000000

In the OEIS, this sequence is defined as ``$a(n)$ is the number of `templates,' or ways of placing a single digit within an $n^2$ by $n^2$ Sudoku puzzle so that all rows, columns, and $n$ by $n$ blocks have exactly one copy of the digit.''

This connection between preference profiles and Sudokus motivated this paper.

\subsection{Sudoku}

A Sudoku grid is a 9 by 9 square that is divided into three vertical \textit{stacks}, and 3 horizontal \textit{bands}. A 3 by 3 square at the intersection of a stack and a band is called a \textit{region}, a \textit{block} or a \textit{box}. A small 1 by 1 square is called a \textit{cell}.

We can expand these definitions to $n^2$ by $n^2$ Sudokus, where the standard Sudoku corresponds to the case $n=3$. The number of completed Sudokus of size $n^2$ by $n^2$ is sequence A107739 in the OEIS \cite{OEIS} database, which starts as follows with the first term corresponding to $n = 0$:
\[1,\ 1,\ 288,\ 6670903752021072936960, \ldots.\]

\subsection{Soulmates, hell-pairs, and outcasts}

The \textit{egalitarian cost} of a man-woman pair is defined as the sum of the rankings they give each other \cite{GI1989}. The concept of egalitarian cost is extended to stable matchings: the egalitarian cost of a stable matching is the sum of the egalitarian costs of all the married couples in the matching. The stable matching with the best possible egalitarian cost for a given preference profile is called \textit{egalitarian}.

By definition, the egalitarian cost of every pair of people in a pseudo-Latin profile is $n + 1$. Therefore, the egalitarian cost of any pseudo-Latin stable matching is $n(n+1)$.

Note that the minimum possible egalitarian cost for a pair of people is $1 + 1 = 2$, because at best, two people can rank each other first. In this scenario, unless they are married to each other, they will be a blocking pair. For this reason, we call such a pair a \textit{soulmate pair}, or simply \textit{soulmates}. When a soulmate pair exists, we can marry them to one another and then reduce the situation to the stable marriage problem with $2n - 2$ people. Since the minimum egalitarian cost is $2$, the egalitarian cost of an egalitarian matching (for $n$ men and $n$ women) is at least $2n$. The cost $2n$ is achievable when all people can be divided into soulmate pairs.

Now, suppose we have two people that rank each other as $n$: the worst they can. We call such a pair a \textit{hell-pair}. If they are married in a stable matching, then they are called a \textit{hell-couple}. Such a pair has an egalitarian cost of $2n$. It follows that such a pair can't happen in a pseudo-Latin preference profile. Also, a stable matching can't contain two hell-couples because two people of opposite genders in different hell-couples form a blocking pair.

In some respect, the opposite notion to soulmates is not a hell-pair, but rather the outcasts. We call two people \textit{outcasts} if they are ranked $n$ by everyone else. In a stable matching, two outcasts have to be matched with each other. Note that the outcasts might be soulmates, hell-couples, or anything in between.

\subsection{The ranking matrices}

We can compress the information about a profile into one $n$ by $n$ matrix, with entry $a_{i,j}$ being $(s,t)$, where $s$ is the ranking by woman $i$ of the man $j$, and $t$ is the ranking by the man $j$ of woman $i$. We call this matrix the \textit{ranking matrix}. If we consider the first number $s$ of every pair $(s,t)$, we get the \textit{women's ranking matrix}. In particular, the first numbers of the pair $(s,t)$ are all different in each row of a ranking matrix. If we instead consider the second number of every pair and transpose the resulting matrix, we get the \textit{men's ranking matrix}. In particular, the second values of the pair $(s,t)$ are all different in each column of a ranking matrix. 

We want to generate another matrix, which ignores the indices of people and only looks at possible pairs of rankings. We call it the \textit{ranking tally matrix}. In other words, if a profile has $k$ pairs of people who rank each other as $(i,j)$, then the element $a_{i,j}$ of the ranking tally matrix equals $k$. The total of the digits in a ranking tally matrix is always $n^2$.

\section{Translating Stable Marriages to Sudoku}\label{sec:SMP2Sudokus}

\subsection{Preference profiles}

To convert preference profiles into a Sudoku grid, we must first split the Sudoku grid into bands representing people of one gender and stacks representing people of the other gender. Without loss of generality, we assume that the stacks represent men and the bands represent women. We number rows in a band from top to bottom and columns in a stack from left to right.

Next, to fill in the cells, we consider an $n$ by $n$ box formed by the intersection of a stack and a band representing a pair of people, say man $M$ and woman $W$.  An entry in the box represents how $M$ and $W$ rank each other; the entry in row $i$ and column $j$ of the box means that woman $W$ ranks man $M$ as $i$ and man $M$ ranks woman $W$ as $j$.

A preference profile marks one cell in each box, and corresponds to how each man-woman pair rank each other.  In our diagrams, we use blue circles to mark such cells. There is exactly one cell marked in each row or column. Indeed, a column corresponds to man $M$ and rank number. There is exactly one person with a given rank by a given man. A similar argument works for a row. This means a preference profile corresponds to one digit filled in a Sudoku grid.

\subsection{Ranking matrices}

A ranking matrix describes the places where blue circles are placed. A matrix entry $a_{i,j} = (s,t)$ says that, in the box that is the intersection of the $i$-th band and $j$-th stack, the blue circle should be in the $s$-th row and $t$-th column. The structure of the matrix matches the structure of Sudoku boxes.

Now we describe the ranking tally matrix in terms of Sudoku. We can imagine that the Sudoku boxes are placed on top of each other in one box. Then, the ranking tally matrix counts the total number of blue circles in each cell.

\subsection{A 2 by 2 example}

\begin{figure}[htp]
    \centering
    \includegraphics[scale=0.1]{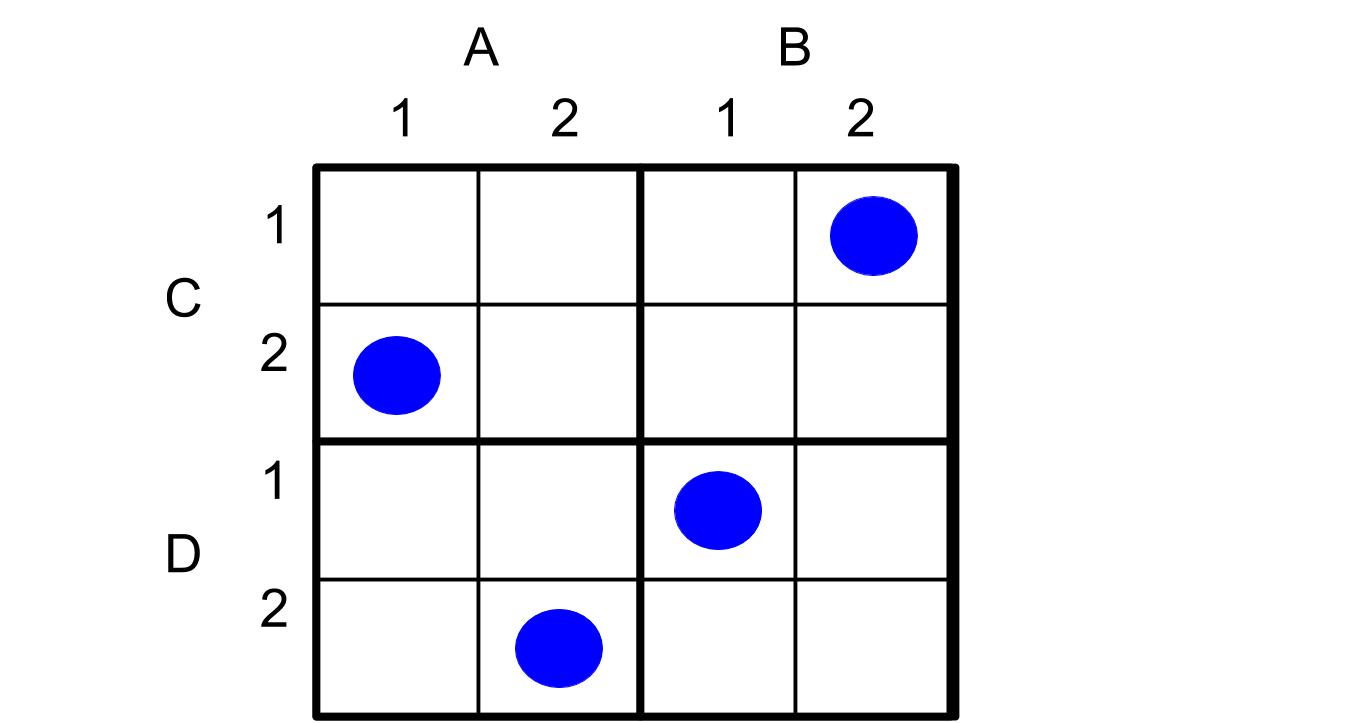}
    \caption{Sudoku example for 2 men and 2 women.}
    \label{fig:SampleProfile2grid}
\end{figure}

We say that $A$ and $B$ are men, and $C$ and $D$ are women. Suppose for example, our preference profile looks like this: $A$ prefers $C$ over $D$, $B$ prefers $D$ over $C$, and both $C$ and $D$ prefer $B$ over $A$. Figure~\ref{fig:SampleProfile2grid} demonstrates how the grid is filled in.

In this example, $B$ and $D$ are soulmates, and in every stable matching, they end up together. It follows that $A$ and $C$ have to be a couple in every stable matching. Thus, this profile allows only one stable matching.

The matrix of how women rank men is in this case $\big(\begin{smallmatrix}
  2 & 1\\
  2 & 1
\end{smallmatrix}\big)$. The matrix of how men rank women is $\big(\begin{smallmatrix}
  1 & 2\\
  2 & 1
\end{smallmatrix}\big)$. 
Here we show the ranking matrix and the ranking tally matrix for this example:
\[
\begin{pmatrix}
  (2,1) & (1,2)\\ 
  (2,2) & (1,1)
\end{pmatrix}
\quad \textrm{ and } \quad
\begin{pmatrix}
  1 & 1\\ 
  1 & 1
\end{pmatrix}.
\]
The tally matrix with all ones is special. Such profiles are called disjoint profiles, and we study them in Section~\ref{sec:DG}.

\subsection{A blocking pair: love gone wrong}

\begin{figure}[htp]
    \centering
    \includegraphics[scale=0.15]{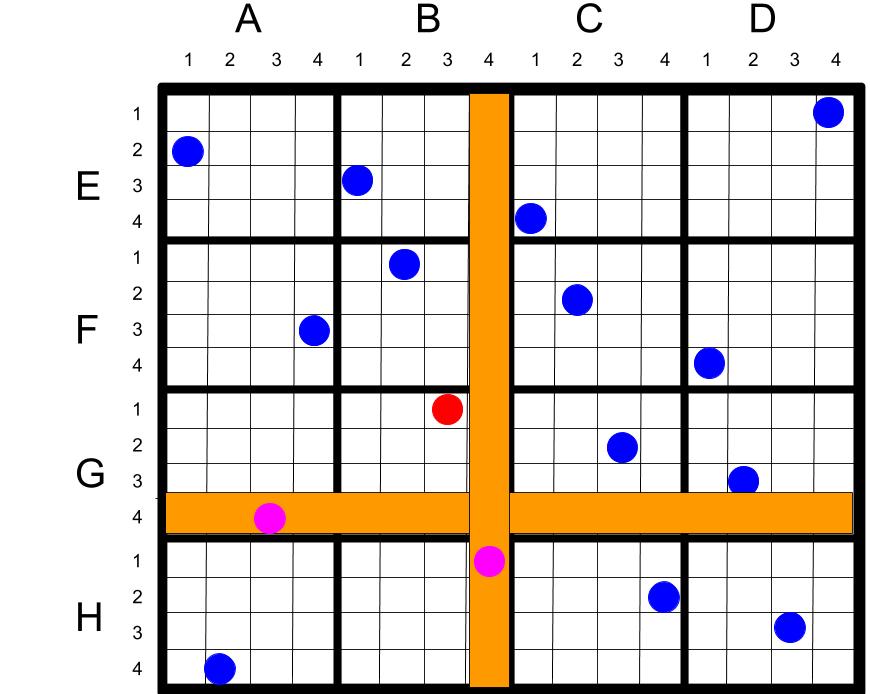}
    \caption{A blocking pair example.}
    \label{fig:square9}
\end{figure}

Figure~\ref{fig:square9} illustrates the notion of a blocking pair. Consider man B and woman G. The pink circles represent their current marriages.  Next, we highlight the column representing B's ranking of his wife and the row representing G's ranking of her husband.  Now, consider the BG box and check the circle corresponding to the mutual ranking of B and G, which is marked red in the figure.  If it is located above and to the left of the shadings, then B and G are a blocking pair. The reasoning is simple. The red circle is to the left of the shaded column, so man B prefers woman G over his current partner since it would be a higher preference. Similarly, the red circle is above the shaded row, so woman G prefers man B over her current partner by the same logic.

To summarize, the blocking pair corresponds to the circle in the Sudoku which is to the left of the circle for the man's current partner and above the circle for the woman's current partner.

\subsection{Gale-Shapley algorithm}

\begin{figure}[ht!]
\begin{center}
\begin{tabular}{ccc}
\centering
    \centering
    \begin{tikzpicture}[scale = 0.5]
        \filldraw[blue] (2.8125, 2.8125) circle (3.5 pt);
        \filldraw[blue] (-2.8125, 2.4375) circle (3.5 pt);
        \filldraw[blue] (-1.3125, 2.0625) circle (3.5 pt);
        \filldraw[blue] (0.1875, 1.6875) circle (3.5 pt);
        \filldraw[blue] (-0.9375, 1.3125) circle (3.5 pt);
        \filldraw[blue] (0.5625, 0.9375) circle (3.5 pt);
        \filldraw[blue] (-1.6875, 0.5625) circle (3.5 pt);
        \filldraw[blue] (1.6875, 0.1875) circle (3.5 pt);
        \filldraw[blue] (-0.5625, -0.1875) circle (3.5 pt);
        \filldraw[blue] (0.9375, -0.5625) circle (3.5 pt);
        \filldraw[blue] (2.0625, -0.9375) circle (3.5 pt);
        \filldraw[blue] (-2.0265, -1.3125) circle (3.5 pt);
        \filldraw[blue] (-0.1875, -1.6875) circle (3.5 pt);
        \filldraw[blue] (1.3125, -2.0635) circle (3.5 pt);
        \filldraw[blue] (2.4375, -2.4375) circle (3.5 pt);
        \filldraw[blue] (-2.4375, -2.8125) circle (3.5 pt);
        \node at (0,-3.5) [text width = 3 cm] {\small Preference Profile};
        \PrefSquareFourSmallEditable
    \end{tikzpicture} & 
    \begin{tikzpicture}[scale = 0.5]
        \filldraw[blue] (2.8125, 2.8125) circle (3.5 pt);
        \filldraw[green] (-2.8125, 2.4375) circle (3.5 pt);
        \filldraw[green] (-1.3125, 2.0625) circle (3.5 pt);
        \filldraw[green] (0.1875, 1.6875) circle (3.5 pt);
        \filldraw[blue] (-0.9375, 1.3125) circle (3.5 pt);
        \filldraw[blue] (0.5625, 0.9375) circle (3.5 pt);
        \filldraw[blue] (-1.6875, 0.5625) circle (3.5 pt);
        \filldraw[green] (1.6875, 0.1875) circle (3.5 pt);
        \filldraw[blue] (-0.5625, -0.1875) circle (3.5 pt);
        \filldraw[blue] (0.9375, -0.5625) circle (3.5 pt);
        \filldraw[blue] (2.0625, -0.9375) circle (3.5 pt);
        \filldraw[blue] (-2.0265, -1.3125) circle (3.5 pt);
        \filldraw[blue] (-0.1875, -1.6875) circle (3.5 pt);
        \filldraw[blue] (1.3125, -2.0635) circle (3.5 pt);
        \filldraw[blue] (2.4375, -2.4375) circle (3.5 pt);
        \filldraw[blue] (-2.4375, -2.8125) circle (3.5 pt);
        \node at (1.75,-3.5) [text width = 3 cm] {\small Round 1};
        \PrefSquareFourSmallEditable
    \end{tikzpicture}& 
    \begin{tikzpicture}[scale = 0.5]
        \filldraw[blue] (2.8125, 2.8125) circle (3.5 pt);
        \filldraw[red] (-2.8125, 2.4375) circle (3.5 pt);
        \filldraw[blue] (-1.3125, 2.0625) circle (3.5 pt);
        \filldraw[blue] (0.1875, 1.6875) circle (3.5 pt);
        \filldraw[blue] (-0.9375, 1.3125) circle (3.5 pt);
        \filldraw[blue] (0.5625, 0.9375) circle (3.5 pt);
        \filldraw[blue] (-1.6875, 0.5625) circle (3.5 pt);
        \filldraw[red] (1.6875, 0.1875) circle (3.5 pt);
        \filldraw[blue] (-0.5625, -0.1875) circle (3.5 pt);
        \filldraw[blue] (0.9375, -0.5625) circle (3.5 pt);
        \filldraw[blue] (2.0625, -0.9375) circle (3.5 pt);
        \filldraw[blue] (-2.0265, -1.3125) circle (3.5 pt);
        \filldraw[blue] (-0.1875, -1.6875) circle (3.5 pt);
        \filldraw[blue] (1.3125, -2.0635) circle (3.5 pt);
        \filldraw[blue] (2.4375, -2.4375) circle (3.5 pt);
        \filldraw[blue] (-2.4375, -2.8125) circle (3.5 pt);
        \node at (0.25,-3.5) [text width = 3 cm] {\small End of Round 1};
        \PrefSquareFourSmallEditable
    \end{tikzpicture}\\
    \begin{tikzpicture}[scale = 0.5]
        \filldraw[blue] (2.8125, 2.8125) circle (3.5 pt);
        \filldraw[red] (-2.8125, 2.4375) circle (3.5 pt);
        \filldraw[blue] (-1.3125, 2.0625) circle (3.5 pt);
        \filldraw[blue] (0.1875, 1.6875) circle (3.5 pt);
        \filldraw[green] (-0.9375, 1.3125) circle (3.5 pt);
        \filldraw[green] (0.5625, 0.9375) circle (3.5 pt);
        \filldraw[blue] (-1.6875, 0.5625) circle (3.5 pt);
        \filldraw[red] (1.6875, 0.1875) circle (3.5 pt);
        \filldraw[blue] (-0.5625, -0.1875) circle (3.5 pt);
        \filldraw[blue] (0.9375, -0.5625) circle (3.5 pt);
        \filldraw[blue] (2.0625, -0.9375) circle (3.5 pt);
        \filldraw[blue] (-2.0265, -1.3125) circle (3.5 pt);
        \filldraw[blue] (-0.1875, -1.6875) circle (3.5 pt);
        \filldraw[blue] (1.3125, -2.0635) circle (3.5 pt);
        \filldraw[blue] (2.4375, -2.4375) circle (3.5 pt);
        \filldraw[blue] (-2.4375, -2.8125) circle (3.5 pt);
        \node at (1.75,-3.5) [text width = 3 cm] {\small Round 2};
        \PrefSquareFourSmallEditable
    \end{tikzpicture} & 
    \begin{tikzpicture}[scale = 0.5]
        \filldraw[blue] (2.8125, 2.8125) circle (3.5 pt);
        \filldraw[red] (-2.8125, 2.4375) circle (3.5 pt);
        \filldraw[blue] (-1.3125, 2.0625) circle (3.5 pt);
        \filldraw[blue] (0.1875, 1.6875) circle (3.5 pt);
        \filldraw[red] (-0.9375, 1.3125) circle (3.5 pt);
        \filldraw[blue] (0.5625, 0.9375) circle (3.5 pt);
        \filldraw[blue] (-1.6875, 0.5625) circle (3.5 pt);
        \filldraw[blue] (1.6875, 0.1875) circle (3.5 pt);
        \filldraw[blue] (-0.5625, -0.1875) circle (3.5 pt);
        \filldraw[blue] (0.9375, -0.5625) circle (3.5 pt);
        \filldraw[blue] (2.0625, -0.9375) circle (3.5 pt);
        \filldraw[blue] (-2.0265, -1.3125) circle (3.5 pt);
        \filldraw[blue] (-0.1875, -1.6875) circle (3.5 pt);
        \filldraw[blue] (1.3125, -2.0635) circle (3.5 pt);
        \filldraw[blue] (2.4375, -2.4375) circle (3.5 pt);
        \filldraw[blue] (-2.4375, -2.8125) circle (3.5 pt);
        \node at (0.25,-3.5) [text width = 3 cm] {\small End of Round 2};
        \PrefSquareFourSmallEditable
    \end{tikzpicture}& 
    \begin{tikzpicture}[scale = 0.5]
        \filldraw[blue] (2.8125, 2.8125) circle (3.5 pt);
        \filldraw[red] (-2.8125, 2.4375) circle (3.5 pt);
        \filldraw[blue] (-1.3125, 2.0625) circle (3.5 pt);
        \filldraw[blue] (0.1875, 1.6875) circle (3.5 pt);
        \filldraw[red] (-0.9375, 1.3125) circle (3.5 pt);
        \filldraw[blue] (0.5625, 0.9375) circle (3.5 pt);
        \filldraw[blue] (-1.6875, 0.5625) circle (3.5 pt);
        \filldraw[blue] (1.6875, 0.1875) circle (3.5 pt);
        \filldraw[blue] (-0.5625, -0.1875) circle (3.5 pt);
        \filldraw[green] (0.9375, -0.5625) circle (3.5 pt);
        \filldraw[green] (2.0625, -0.9375) circle (3.5 pt);
        \filldraw[blue] (-2.0265, -1.3125) circle (3.5 pt);
        \filldraw[blue] (-0.1875, -1.6875) circle (3.5 pt);
        \filldraw[blue] (1.3125, -2.0635) circle (3.5 pt);
        \filldraw[blue] (2.4375, -2.4375) circle (3.5 pt);
        \filldraw[blue] (-2.4375, -2.8125) circle (3.5 pt);
        \node at (1.75,-3.5) [text width = 3 cm] {\small Round 3};
        \PrefSquareFourSmallEditable
    \end{tikzpicture}\\
    \begin{tikzpicture}[scale = 0.5]
        \filldraw[blue] (2.8125, 2.8125) circle (3.5 pt);
        \filldraw[red] (-2.8125, 2.4375) circle (3.5 pt);
        \filldraw[blue] (-1.3125, 2.0625) circle (3.5 pt);
        \filldraw[blue] (0.1875, 1.6875) circle (3.5 pt);
        \filldraw[red] (-0.9375, 1.3125) circle (3.5 pt);
        \filldraw[blue] (0.5625, 0.9375) circle (3.5 pt);
        \filldraw[blue] (-1.6875, 0.5625) circle (3.5 pt);
        \filldraw[blue] (1.6875, 0.1875) circle (3.5 pt);
        \filldraw[blue] (-0.5625, -0.1875) circle (3.5 pt);
        \filldraw[red] (0.9375, -0.5625) circle (3.5 pt);
        \filldraw[blue] (2.0625, -0.9375) circle (3.5 pt);
        \filldraw[blue] (-2.0265, -1.3125) circle (3.5 pt);
        \filldraw[blue] (-0.1875, -1.6875) circle (3.5 pt);
        \filldraw[blue] (1.3125, -2.0635) circle (3.5 pt);
        \filldraw[blue] (2.4375, -2.4375) circle (3.5 pt);
        \filldraw[blue] (-2.4375, -2.8125) circle (3.5 pt);
        \node at (0.25,-3.5) [text width = 3 cm] {\small End of Round 3};
        \PrefSquareFourSmallEditable
    \end{tikzpicture}&
    \begin{tikzpicture}[scale = 0.5]
        \filldraw[blue] (2.8125, 2.8125) circle (3.5 pt);
        \filldraw[red] (-2.8125, 2.4375) circle (3.5 pt);
        \filldraw[blue] (-1.3125, 2.0625) circle (3.5 pt);
        \filldraw[blue] (0.1875, 1.6875) circle (3.5 pt);
        \filldraw[red] (-0.9375, 1.3125) circle (3.5 pt);
        \filldraw[blue] (0.5625, 0.9375) circle (3.5 pt);
        \filldraw[blue] (-1.6875, 0.5625) circle (3.5 pt);
        \filldraw[blue] (1.6875, 0.1875) circle (3.5 pt);
        \filldraw[blue] (-0.5625, -0.1875) circle (3.5 pt);
        \filldraw[red] (0.9375, -0.5625) circle (3.5 pt);
        \filldraw[blue] (2.0625, -0.9375) circle (3.5 pt);
        \filldraw[blue] (-2.0265, -1.3125) circle (3.5 pt);
        \filldraw[blue] (-0.1875, -1.6875) circle (3.5 pt);
        \filldraw[blue] (1.3125, -2.0635) circle (3.5 pt);
        \filldraw[green] (2.4375, -2.4375) circle (3.5 pt);
        \filldraw[blue] (-2.4375, -2.8125) circle (3.5 pt);
        \node at (1.75,-3.5) [text width = 3 cm] {\small Round 4};
        \PrefSquareFourSmallEditable
    \end{tikzpicture}&
    \begin{tikzpicture}[scale = 0.5]
        \filldraw[blue] (2.8125, 2.8125) circle (3.5 pt);
        \filldraw[red] (-2.8125, 2.4375) circle (3.5 pt);
        \filldraw[blue] (-1.3125, 2.0625) circle (3.5 pt);
        \filldraw[blue] (0.1875, 1.6875) circle (3.5 pt);
        \filldraw[red] (-0.9375, 1.3125) circle (3.5 pt);
        \filldraw[blue] (0.5625, 0.9375) circle (3.5 pt);
        \filldraw[blue] (-1.6875, 0.5625) circle (3.5 pt);
        \filldraw[blue] (1.6875, 0.1875) circle (3.5 pt);
        \filldraw[blue] (-0.5625, -0.1875) circle (3.5 pt);
        \filldraw[red] (0.9375, -0.5625) circle (3.5 pt);
        \filldraw[blue] (2.0625, -0.9375) circle (3.5 pt);
        \filldraw[blue] (-2.0265, -1.3125) circle (3.5 pt);
        \filldraw[blue] (-0.1875, -1.6875) circle (3.5 pt);
        \filldraw[blue] (1.3125, -2.0635) circle (3.5 pt);
        \filldraw[red] (2.4375, -2.4375) circle (3.5 pt);
        \filldraw[blue] (-2.4375, -2.8125) circle (3.5 pt);
        \node at (0.25,-3.5) [text width = 3 cm] {\small End of Round 4};
        \PrefSquareFourSmallEditable
    \end{tikzpicture}
\end{tabular}
\end{center}
    \caption{Gale-Shapley algorithm example.}
    \label{fig:GSae}
\end{figure}
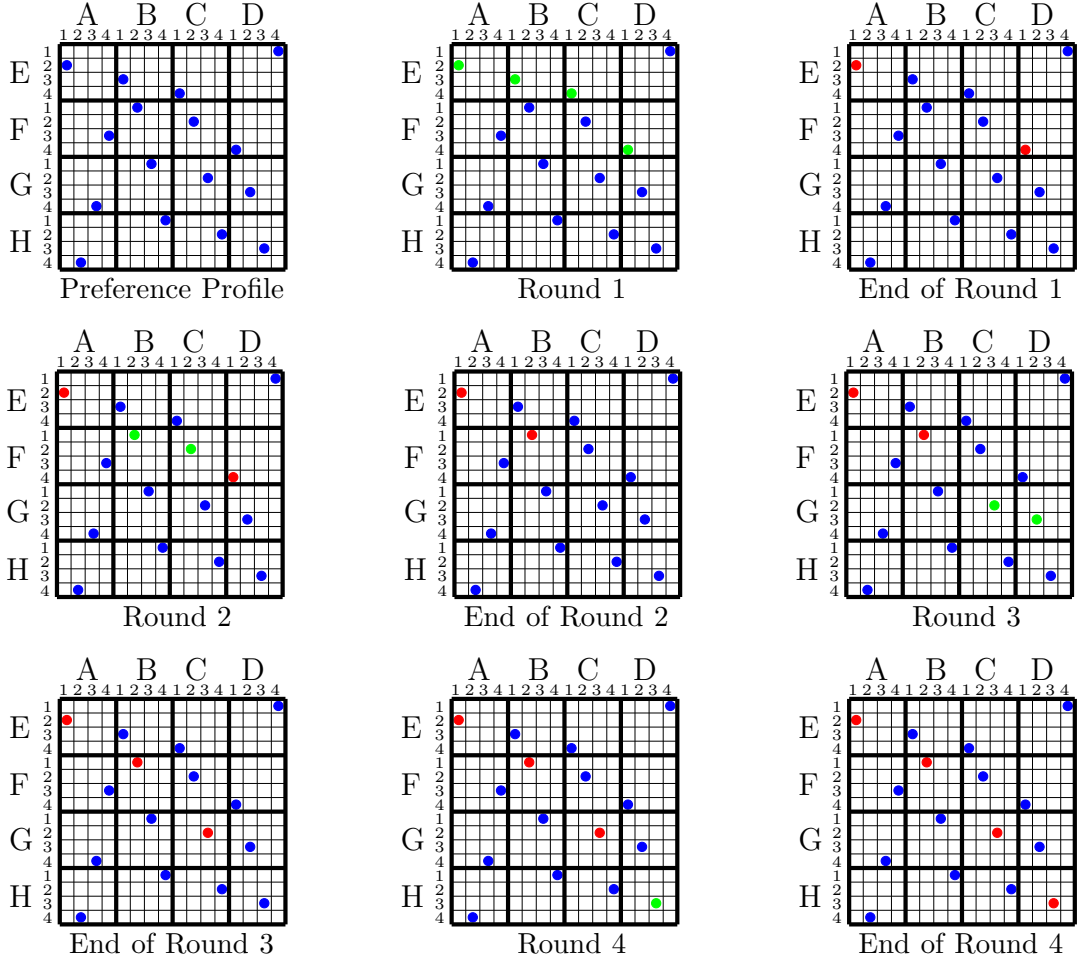

We illustrate the men-proposing Gale-Shapley algorithm in terms of Sudoku in Figure~\ref{fig:GSae}. The first round starts with all the circles colored blue. In the process, green circles correspond to proposals, and red circles correspond to engagements. At the end of each round, there are no green circles.

In the first round, men propose to the women they rank first. This means we color green the circle in the first column of every stack. Now women choose the best proposals. If there are several green circles in a band, then we color the entry in the top-most row red and revert other circles, if any, to blue. The red circles show engagements at the end of the round.

Then we repeat the following steps until every stack has a red circle. If a stack has a red circle, meaning the corresponding man is engaged, then the stack is unchanged. If a stack doesn't have a red circle, then we find the left-most circle that has never been recolored and color it green. In other words, unengaged men propose to their next best choice. 

Now we go to bands (women's turn). The bands without green circles (the bands without new proposals) are unchanged. If a band has a green circle, then we color the top-most non-blue circle red, and reset other circles to blue. In the end, we have one red circle in each band and stack corresponding to a stable matching.

\subsection{The egalitarian cost}

In terms of Sudoku, the egalitarian cost of a marriage between $M$ and $W$ can be defined using the entry in the box corresponding to $M$ and $W$. Namely, it is the Manhattan distance from the entry for their preferences to the top-left corner plus 2. Note that the Manhattan distance between two cells in a grid is the smallest number of steps from one cell to the other. In this case, an allowed step is a move from a cell to an orthogonally-adjacent cell. The egalitarian cost of a cell in a block could also be determined by the anti-diagonal it is in since the anti-diagonal is equidistant from the top-left cell.

The egalitarian cost of a stable matching is $2n$ plus the sum of the Manhattan distances between each red circle (circles correspond to couples in the matching) and the top-left corner of the box.

\section{Translating Sudokus to stable marriages}\label{sec:Sudokus2SMP}

\subsection{Complete Sudokus}

Cells corresponding to a particular digit in a complete Sudoku grid correspond to a preference profile. The profiles corresponding to different digits and the same two people must have a different pair of mutual rankings, as the digits must occupy different cells in a box. In other words, for any two profiles in any man-woman pair, the ranking of the woman in the man's profiles must be different, or the ranking of the man in the woman's profiles must be different. We call two profiles \textit{non-overlapping} if they correspond to two different digits in a complete Sudoku grid. Thus, nine non-overlapping profiles correspond to all nine digits in a complete Sudoku grid.

One Sudoku clue means placement of one digit in a Sudoku grid. It corresponds to a pair of a man and a woman and their ranking of each other. A Sudoku puzzle means that there are given several mutual rankings for several non-overlapping profiles. Finishing a Sudoku puzzle translates into extending given mutual rankings to nine complete non-overlapping profiles. 

The discussion below only relates to $n=3$, aka standard Sudoku.

\subsection{Total number of Sudokus}

Given that the number of different Sudoku grids is 6670903752021072936960 \cite{RT2011}, the number of nine non-overlapping profiles is the same.

\subsection{Minimum number of clues}

The minimum number of clues that a Sudoku puzzle can have and still produce a unique solution is 17 \cite{RT2011}. This number was found by an exhaustive search. Not all Sudoku grids with 17 clues have unique solutions, but no Sudoku grid with fewer than 17 clues has a unique solution. In terms of the Stable Marriage Problem, this means that we need at least 17 non-overlapping mutual ranking values to complete the nine non-overlapping profiles. Still, there is no guarantee that we can determine all 9 non-overlapping profiles even if we do have 17 mutual ranking values.

\subsection{Maximum number of clues}

What is the maximum number of clues that do not guarantee a unique solution? One can see that if the puzzle is missing one entry of a given digit, this digit has a unique placement. If the missing clues are of the same digit, the solution is unique too. By combining these statements, we see that if not more than one digit has more than one entry missing, we can finish the puzzle uniquely. It follows that if there are 80, 79, or 78 clues given, the solution is unique. However, it is possible to construct a Sudoku puzzle with 77 clues given and two solutions.

By extending this example to any $n$, we see that $n^4-3$ clues guarantee a unique solution to a puzzle if the solution exists. Correspondingly, we see that $n^4-3$ rankings are necessary to guarantee that one can reconstruct $n^2$ unique non-overlapping preference profiles.

\section{Special profiles}\label{sec:specialprofiles}

\subsection{Generalized pseudo-Latin profiles}

Recall that pseudo-Latin profiles \cite{T2002} are profiles where the egalitarian cost of each pair of people is $n+1$. In such profiles, two men can't give the same rank to the same woman; otherwise, she would have ranked both men the same. As a consequence, the matrix of men's preferences forms a Latin square. Similarly, the matrix of women's preferences forms a Latin square.

It is natural to drop the request that the egalitarian cost for every pair of people is $n+1$ while leaving the requirement that both the men's profiles and the women's profiles form a Latin square. We call such profiles \textit{mutually Latin} profiles.

The Gale-Shapley algorithm, when applied on a mutually Latin profile, will always take exactly one round. Because every man has a different ranking for each woman, all the men propose to different women in the first round. However, the number of stable matchings can vary from one (when there are $n$ pairs of soulmates) to relatively large when it is a pseudo-Latin profile.

For $n > 2$, we can't have outcasts in a mutually Latin profile, as we have exactly one man who ranks a given woman as $n$ and vice versa. But we can have soulmates and hell-pairs. We can even have $n$ pairs of soulmates and $n$ hell-pairs.

\subsection{Disjoint-groups Sudoku and disjoint profiles}\label{sec:DG}

There is a special type of Sudoku, where in a particular place in a box, all digits are distinct across the boxes. For example, the top-left corners of every box have distinct digits. This creates 9 additional groups to add to columns, rows, and boxes that have to contain distinct digits.  Such Sudokus are called \textit{disjoint-groups Sudokus} or DG Sudokus. This Sudoku type is used a lot in Sudoku puzzles.

A Sudoku corresponds to $n$ non-overlapping profiles. Each profile in a DG Sudoku has a property that each pair of mutual rankings $(i,j)$ occurs exactly once. In other words, the ranking tally matrix consists of all ones. We call such a profile a \textit{disjoint} profile. In particular, each disjoint profile has exactly one pair of soulmates and one hell-pair. 

Figure~\ref{fig:LexiDGexample} shows an example of a disjoint-groups Sudoku. Coincidentally, this is the lexicographically earliest complete Sudoku grid, read by rows. 

\begin{figure}[ht!]
\begin{center}
\begin{tikzpicture}
    \TripleEdit
    \FirstColumnTripleEdit{1}{4}{7}{2}{3}{8}{5}{6}{9}
    \SecondColumnTripleEdit{2}{5}{8}{1}{6}{9}{3}{4}{7}
    \ThirdColumnTripleEdit{3}{6}{9}{4}{5}{7}{1}{2}{8}
    \FourthColumnTripleEdit{4}{7}{1}{3}{8}{2}{6}{9}{5}
    \FifthColumnTripleEdit{5}{8}{2}{6}{9}{1}{4}{7}{3}
    \SixthColumnTripleEdit{6}{9}{3}{5}{7}{4}{2}{8}{1}
    \SeventhColumnTripleEdit{7}{1}{4}{8}{2}{3}{9}{5}{6}
    \EighthColumnTripleEdit{8}{2}{5}{9}{1}{6}{7}{3}{4}
    \NinthColumnTripleEdit{9}{3}{6}{7}{4}{5}{8}{1}{2}
\end{tikzpicture}
\end{center}
    \caption{Lexicographically earliest Sudoku, which coincidentally is also a DG Sudoku.}
    \label{fig:LexiDGexample}
\end{figure}

For another example of a disjoint profile, we consider a profile where all men rank everyone the same way. To prove that such a profile is disjoint, consider the mutual ranking $(i,j)$. There is exactly one woman that is ranked $i$ (by every man). And she ranks precisely one man as $j$. Thus, a pair $(i,j)$ appears exactly once in the ranking. 

We can have a Sudoku consisting of profiles where all men rank everyone the same way. An example of such a profile is shown in Figure~\ref{fig:DGexampleAllTheSame}.

\begin{figure}[ht!]
\begin{center}
\begin{tikzpicture}
    \TripleEdit
    \FirstColumnTripleEdit{1}{4}{7}{2}{5}{8}{3}{6}{9}
    \SecondColumnTripleEdit{2}{5}{8}{3}{6}{9}{4}{7}{1}
    \ThirdColumnTripleEdit{3}{6}{9}{4}{7}{1}{5}{8}{2}
    \FourthColumnTripleEdit{4}{7}{1}{5}{8}{2}{6}{9}{3}
    \FifthColumnTripleEdit{5}{8}{2}{6}{9}{3}{7}{1}{4}
    \SixthColumnTripleEdit{6}{9}{3}{7}{1}{4}{8}{2}{5}
    \SeventhColumnTripleEdit{7}{1}{4}{8}{2}{5}{9}{3}{6}
    \EighthColumnTripleEdit{8}{2}{5}{9}{3}{6}{1}{4}{7}
    \NinthColumnTripleEdit{9}{3}{6}{1}{4}{7}{2}{5}{8}
\end{tikzpicture} 
\end{center}
    \caption{An example of a DG Sudoku, where in each corresponding profile all men rank everyone the same way.}
    \label{fig:DGexampleAllTheSame}
\end{figure}

How can we see in a Sudoku that all men have the same preferences? For simplicity, let's assume we are talking about profile 1. In each band, all the ones must all be in the same column in their box. Now we connect such profiles to stable matchings.

\begin{proposition}
If all the people of the same gender have the same list of preferences, there is only one stable matching.
\end{proposition}

\begin{proof}
Without loss of generality, we assume that men have the same preferences. We continue by induction. If $n=1$, we have only one matching possible, and it is stable. Suppose for $n \leq k$ the statement is true. Consider $n =k+1$. The pair of people that consists of the most desirable woman and her first choice are soulmates. They are married in any stable matching. If we remove them from consideration, we can use induction.
\end{proof}

Consider the lexicographically earliest Sudoku in Figure~\ref{fig:LexiDGexample}. For every digit, the men have the same preferences. It follows that each profile has exactly 1 stable matching. By shuffling the digits, we can get $9!$ DG Sudokus where each profile has exactly 1 stable matching. One of these Sudokus is the lexicographically last Sudoku: it can be written out by replacing digit $x$ with digit $10-x$ in the lexicographically earliest Sudoku.

Each disjoint profile has exactly one pair of soulmates. It means that the number of stable matchings for such a profile can't exceed the maximum number of possible matchings for $n-1$ men and women. For example, when $n=3$, each disjoint profile can't have more than 2 stable matchings. In the examples above (Figures~\ref{fig:LexiDGexample} and~\ref{fig:DGexampleAllTheSame}), we have DG Sudokus such that all of the profiles have exactly one stable matching. Figure~\ref{fig:DGexample2StableMatchings} shows an example of a DG Sudoku where two of the profiles, corresponding to digits 1 and 3, have two stable matchings. The other profiles have one stable matching.

\begin{figure}[ht!]
\begin{center}
\begin{tikzpicture}
    \TripleEdit
    \FirstColumnTripleEdit{1}{2}{3}{4}{5}{6}{7}{8}{9}
    \SecondColumnTripleEdit{4}{5}{6}{7}{8}{9}{3}{1}{2}
    \ThirdColumnTripleEdit{7}{8}{9}{2}{3}{1}{4}{5}{6}
    \FourthColumnTripleEdit{6}{4}{5}{3}{1}{2}{9}{7}{8}
    \FifthColumnTripleEdit{9}{7}{8}{6}{4}{5}{1}{2}{3}
    \SixthColumnTripleEdit{3}{1}{2}{9}{7}{8}{6}{4}{5}
    \SeventhColumnTripleEdit{8}{9}{7}{5}{6}{4}{2}{3}{1}
    \EighthColumnTripleEdit{2}{3}{1}{8}{9}{7}{5}{6}{4}
    \NinthColumnTripleEdit{5}{6}{4}{1}{2}{3}{8}{9}{7}
\end{tikzpicture} 
\end{center}
    \caption{An example of a DG Sudoku where some of the profiles have two stable matchings.}
    \label{fig:DGexample2StableMatchings}
\end{figure}

What happens if we combine disjoint profiles with mutually Latin profiles? In mutually Latin profiles men's and women's preferences form Latin squares. By definition of disjoint profiles, when these squares are superimposed, the ordered paired entries in the cells of the square are all distinct. This is exactly the definition of the \textit{mutually orthogonal} or \textit{Graeco-Latin} squares.

Mutually-orthogonal squares are used to build magic squares. Thus, we can use the men's and women's preferences in a disjoint mutually Latin profile to build a magic square.

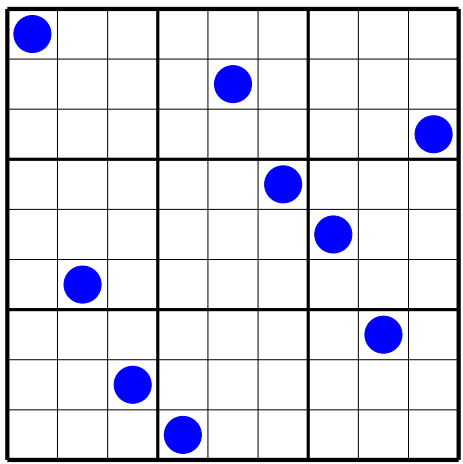
\begin{figure}[ht!]
\begin{center}
\begin{tikzpicture}
    \TripleEdit
    \filldraw[blue] (-8/3,8/3) circle (7 pt);
    \filldraw[blue] (0,6/3) circle (7 pt);
    \filldraw[blue] (8/3,4/3) circle (7 pt);
    \filldraw[blue] (-6/3,-2/3) circle (7 pt);
    \filldraw[blue] (2/3,2/3) circle (7 pt);
    \filldraw[blue] (4/3,0/3) circle (7 pt);
    \filldraw[blue] (-4/3,-6/3) circle (7 pt);
    \filldraw[blue] (-2/3,-8/3) circle (7 pt);
    \filldraw[blue] (6/3,-4/3) circle (7 pt);
\end{tikzpicture}
\end{center}
    \caption{A disjoint mutually Latin profile.}
    \label{fig:pLD}
\end{figure}

Figure~\ref{fig:pLD} shows an example of a disjoint mutually Latin profile. Here are the preference matrices for men and women in this profile:
\[\begin{pmatrix}
 1  & 2 & 3\\ 
 3  & 1 & 2 \\
 2  & 3 & 1
\end{pmatrix}
\quad \textrm{ and } \quad
\begin{pmatrix}
 1  & 2 & 3\\ 
 2  & 3 & 1 \\
 3  & 1 & 2
\end{pmatrix}.
\]
The first man and the first woman are soulmates, so they have to be a couple in any stable matching. For the remaining four people, both possible matchings are stable. Therefore, the total number of stable matching for this profile is 2. 

We gave an example above (where all men have the same preferences) of a disjoint profile that is not a mutually Latin profile. An example of a mutually Latin profile that is not disjoint exists. The following subsection discusses such an example.

\section{Joint-groups Sudoku and joint profiles}\label{sec:joint}

We introduce a new type of profile that is, in a sense, the opposite of a disjoint profile. In a disjoint profile, every pair of rankings is possible and occurs once. In this new type of profile, we want every pair of rankings to repeat as much as possible. Since a pair of rankings can appear in a profile not more than $n$ times, in this new profile type, each possible ranking $(i,j)$ appears exactly $n$ times, and there is one value $j$ for every $i$. For example, if a number appears in the top-left corner of a box, it has to appear in the top-left corner of $n$ boxes. We call such a profile a \textit{joint} profile. We call a complete Sudoku with each digit forming a joint profile a \textit{joint-groups} Sudoku.

\subsection{Joint profiles}

The ranking tally matrix of such a joint profile has $n$ non-zero elements, and each of them equal to $n$. Non-zero elements form an entry in a Latin square: that is, there is exactly one non-zero value in each column and each row. 

Suppose $(i,j)$ is a pair of mutual rankings between a man and a woman. In a joint profile, $j$ is a function of $i$: $j = f(i)$. That means if a man is ranked $i$ by a woman, he has to rank her back as $f(i)$. We call the function $f$ the \textit{key}.

For example, consider a special profile, where $f(i) = i$. We call it a \textit{mirror} profile. The corresponding ranking tally matrix is a diagonal matrix: it is an identity matrix multiplied by $n$. A mirror profile has $n$ pairs of soulmates and one stable matching.

A joint profile is a relaxation of a pseudo-Latin profile. Indeed, the key function for a pseudo-Latin profile is $f(i) = n+1-i$.

\begin{proposition}
The key $f(i)$ is a bijection.
\end{proposition}

\begin{proof}
Suppose $f(i)$ is not a bijection. That means, there are two distinct values $i_1$ and $i_2$ such that $f(i_1) = f(i_2) = j$. It follows that men are ranked $j$ at least $2n$ times. Thus there exists a man that is ranked $j$ by at least two women. Hence, he ranks both of these women the same as $f(j)$, which contradicts the notion of ranking. 
\end{proof}

\begin{theorem}
A joint profile is a mutually Latin profile.
\end{theorem}

\begin{proof}
Because $f(i)$ is a bijection, two women cannot rank a particular man the same way. Thus the women's profiles form a Latin square. Moreover, as $f(i)$ is a bijection, its inverse is the key function for men. Thus, the men's profiles form a Latin square too.
\end{proof}

Suppose $n=3$. Up to relabeling the men, the first woman's preferences are 1,2,3. Then, we can swap the other two women to guarantee that the second woman's preferences are 2,3,1 and the third woman's preferences are 3,1,2. That means that, up to symmetries, all women's preferences are equivalent to each other. Hence, the function $f$, up to symmetries, defines the profile uniquely and, consequently, uniquely determines the number of stable matchings.

Joint profiles are mutually Latin profiles, but not vice versa. Figure~\ref{fig:pLD} shows a mutually Latin disjoint profile. Thus, there exists a mutually Latin profile that is not a joint profile. 

Unlike the disjoint profiles, where the total number of stable matchings is less than the potential maximum number of possible stable matchings, the number of stable matchings in a joint profile can vary greatly. If $f(i) = i$, we get a mirror profile with $n$ pairs of soulmates and exactly one stable matching. On the other hand, if $f(i) = n+1-i$, we get a pseudo-Latin profile, and they tend to produce a large number of stable matchings.

Now we want to move our attention to hell-pairs in joint profiles.

\begin{proposition}
No hell-couples can exist in a stable matching for a joint profile with $n$ men and $n$ women.
\end{proposition}

\begin{proof}
Suppose there exists a stable matching containing a hell-couple for a given joint profile. That means $f(n) = n$, where $f$ is the key function. Suppose $M$ is the man in the couple. Because the women's preferences form a Latin square, there exists a woman $W$ who ranks man $M$ first. Then $M$ and $W$ form a blocking pair. This leads to a contradiction.
\end{proof}

\subsection{Uniform profiles}

Suppose we have a joint profile, and the mutual ranking $(i,j)$ is possible. Then we can match all people into pairs with this mutual ranking. We call such a matching \textit{uniform}. The following proposition describes when a uniform matching is stable.

\begin{proposition}\label{prop:uniform}
Given a joint profile, a uniform matching corresponding to mutual ranking $(i,j)$ is unstable if there exists another mutual ranking $(k,\ell)$ such that $k < i$ and $\ell < j$. Otherwise, it is stable.
\end{proposition}

\begin{proof}
Consider a uniform matching corresponding to mutual ranking $(i,j)$. Such a matching is unstable if and only if a blocking pair exists, that is, a pair with mutual ranking $(k,\ell)$ such that $k < i$ and $\ell < j$. Suppose such ranking $(k,\ell)$ exists. Given that our matching is uniform, any pair of people with mutual ranking $(k,\ell)$ are not matched to each other and thus form a blocking pair.
\end{proof}

For example, the following uniform matchings are stable in a joint profile:
\begin{itemize}
\item All couples with ranking $(1,k)$. In this case, the women-proposing Gale-Shapley algorithm ends in one round. It is a stable matching that is woman-optimal.
\item All couples with ranking $(j,1)$. In this case, the men-proposing Gale-Shapley algorithm ends in one round. It is a stable matching that is man-optimal.
\item All couples with ranking $(i,j)$, where $i+j$ is the lowest egalitarian cost. This is an egalitarian matching.
\end{itemize}

\begin{corollary}
For a pseudo-Latin profile, there are at least $n$ different stable matchings.
\end{corollary}

\begin{proof}
Each possible mutual ranking has the same egalitarian cost. Thus, all pairs with the same mutual ranking form a stable matching.
\end{proof}

When $n$ increases, it is possible to have a profile such that not all couples have the same mutual ranking for each other. An example of such a profile is shown below, represented as a ranking matrix:
\[
\begin{pmatrix}
  (3,4) & (4,3) & \textbf{(1,2)} & (2,1)\\ 
  (4,3) & (3,4) & (2,1) & \textbf{(1,2)}\\
  (1,2) & \textbf{(2,1)} & (3,4) & (4,3) \\
  \textbf{(2,1)} & (1,2) & (4,3) & (3,4)
\end{pmatrix}.
\]
A stable matching is in bold. It has two types of mutual rankings: (1,2) and (2,1).

\subsection{Classification of joint profiles for $n=3$}

\begin{figure}[ht!]
\begin{center}
\begin{tabular}{cccccc}
\centering
    \centering
        \begin{tikzpicture}[scale = 0.3]
            \TripleSquareEdit
            \filldraw[green] (-2,2) circle (22.5pt);
            \filldraw[red] (0,0) circle (22.5pt);
            \filldraw[red] (2,-2) circle (22.5pt);
        \end{tikzpicture} &
        \begin{tikzpicture}[scale = 0.3]
            \TripleSquareEdit
            \filldraw[green] (-2,2) circle (22.5pt);
            \filldraw[red] (2,0) circle (22.5pt);
            \filldraw[red] (0,-2) circle (22.5pt);
        \end{tikzpicture} &
        \begin{tikzpicture}[scale = 0.3]
            \TripleSquareEdit
            \filldraw[green] (2,2) circle (22.5pt);
            \filldraw[green] (0,0) circle (22.5pt);
            \filldraw[green] (-2,-2) circle (22.5pt);
        \end{tikzpicture} &
        \begin{tikzpicture}[scale = 0.3]
            \TripleSquareEdit
            \filldraw[green] (-2,0) circle (22.5pt);
            \filldraw[green] (2,2) circle (22.5pt);
            \filldraw[red] (0,-2) circle (22.5pt);
        \end{tikzpicture} &
        \begin{tikzpicture}[scale = 0.3]
            \TripleSquareEdit
            \filldraw[green] (-2,-2) circle (22.5pt);
            \filldraw[green] (0,2) circle (22.5pt);
            \filldraw[red] (2,0) circle (22.5pt);
        \end{tikzpicture} &
        \begin{tikzpicture}[scale = 0.3]
            \TripleSquareEdit
            \filldraw[green] (-2,0) circle (22.5pt);
            \filldraw[green] (0,2) circle (22.5pt);
            \filldraw[red] (2,-2) circle (22.5pt);
        \end{tikzpicture} \\
\end{tabular}
\end{center}
    \caption{Ranking tally matrices.}
    \label{fig:rankingtally2}
\end{figure}
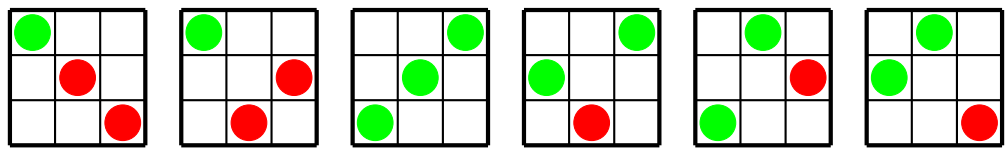

For this section, we assume that $n=3$. We classify joint profiles with respect to the number of stable matchings. We look at the profiles in terms of the key functions. There are six different key functions. The surprising fact is that all stable matchings for $n=3$ are uniform. In other words, they correspond to one entry in the tally matrix. These profiles are represented in Figure~\ref{fig:rankingtally2} by ranking tally matrices, where the digit 3 is replaced with a circle. Green circles show a possible mutual ranking for all couples in a stable matching. By contrast, red circles show the impossible rankings. As we proved in Proposition~\ref{prop:uniform}, the red circles for uniform matchings are the ones that have another circle to the top-left.

Suppose $f(1) = 1$, or in other words, $(1,1)$ belongs to the list of rankings. This means that we have three pairs of soulmates and one stable matching. In this matching, all pairs of soulmates are matched to each other. This case covers two key functions corresponding to the two left-most matrices in Figure~\ref{fig:rankingtally2}.

Suppose we have a pseudo-Latin profile, or $f(i) = 4-i$. Then there are three stable matchings. One matching marries pairs with mutual ranking $(1,3)$: it is woman-optimal and man-pessimal. Similarly, there is a stable matching marrying pairs with mutual ranking $(3,1)$: it is woman-pessimal and man-optimal. The third matching marries pairs with mutual ranking $(2,2)$. This key function corresponds to the third matrix in Figure~\ref{fig:rankingtally2}.

There are three possible key functions left, with two different stable matchings each. One matching marries pairs with mutual ranking $(1,i)$: it is woman-optimal and man-pessimal. Similarly, there is a stable matching marrying pairs with mutual ranking $(i,1)$: it is woman-pessimal and man-optimal. The third matching is impossible in a stable marriage. These three key functions correspond to the last three matrices in Figure~\ref{fig:rankingtally2}.

\subsection{Joint-groups Sudoku}

As we mentioned, we can use a ranking tally matrix to represent a joint profile. It has exactly one non-zero element in a row or column. If two profiles have non-zero elements in different cells of the tally matrix, then they are non-overlapping.

We want to prove that joint-groups Sudokus exist for any $n$. In the theorem below, we build such a Sudoku by cycling boxes. A Sudoku is an $n$ by $n$ grid of boxes. We can consider the first band as a set of $n$ boxes. 

We introduce a \textit{box-cyclic} Sudoku where the boxes for every band are cycled one box to the right compared to the band above. Thus in such a Sudoku, the box in the intersection of band $i$ and stack $j$ is the same as the box in the intersection of band $i+k$ and stack $j+k$, where the band and the stack numbers are taken modulo $n$.

Box-cyclic Sudokus are JG Sudokus. Indeed, in such a Sudoku, each digit appears in exactly $n$ places inside a box. In addition, box-cyclic Sudokus are easy to build. We use this fact in the next theorem.

\begin{theorem}
A joint-groups Sudoku exists for any $n$.
\end{theorem}

\begin{proof}
First, we arrange the integers 1 through $n^2$ into an $n$ by $n$ box. We call it box $B_1$. For box $B_i$, we take the numbers from box $B_1$ and move them $i$ steps to the right and $i$ steps down, wrapping around if necessary.  Now we place box $B_i$ into the intersection of stack $j$ and band $j+i-1$ considered modulo $n$. Thus, every box appears $n$ times: once in each stack and  once in each band. 

Now we show that we get a complete Sudoku grid. Consider a number that appears at coordinates $(a,b)$ in box $B_1$. Then it appears at coordinates $(a+i, b+i)$ modulo $n$ in box $B_i$. Each box appears once in every stack and band; therefore, the number appears once in every row and column, as well as once in every box.

Now we show that the rules of joint-groups Sudoku are satisfied. The key function corresponding to the integer that appears at coordinates $(a,b)$ in box $B_1$ is the following: $f(x) = b-a +x$ modulo $n$.
\end{proof}

A JG Sudoku doesn't have to be box-cyclic, as demonstrated in Figure~\ref{fig:JGbutNotbc}.
\begin{figure}[ht!]
\begin{center}
\begin{tikzpicture}
     \TripleEdit
    \FirstColumnTripleEdit{1}{8}{6}{3}{9}{4}{2}{7}{5}
    \SecondColumnTripleEdit{4}{2}{9}{5}{1}{7}{6}{3}{8}
    \ThirdColumnTripleEdit{7}{5}{3}{8}{6}{2}{9}{4}{1}
    \FourthColumnTripleEdit{2}{9}{4}{1}{7}{5}{3}{8}{6}
    \FifthColumnTripleEdit{5}{3}{7}{6}{2}{8}{4}{1}{9}
    \SixthColumnTripleEdit{8}{6}{1}{9}{4}{3}{7}{5}{2}
    \SeventhColumnTripleEdit{3}{7}{5}{2}{8}{6}{1}{9}{4}
    \EighthColumnTripleEdit{6}{1}{8}{4}{3}{9}{5}{2}{7}
    \NinthColumnTripleEdit{9}{4}{2}{7}{5}{1}{8}{6}{3}
\end{tikzpicture}
\end{center}
    \caption{A JG but not box-cyclic Sudoku.}
    \label{fig:JGbutNotbc}
\end{figure}

Now we want to define the placement matrix for a box-cyclic Sudoku. The \textit{placement matrix} is an $n$ by $n$ matrix, where each element of the matrix shows $n$ integers that appear in the corresponding cells of the box. For example, the placement matrix for Sudoku in Figure~\ref{fig:JGbutNotbc} is
\[\begin{pmatrix}
123 & 456 & 789\\ 
789 & 123 & 456\\
456 & 789 & 123
\end{pmatrix}.
\]

\subsection{JG Sudoku puzzles}

We invented a new type of Sudoku, so we have no choice but to present some puzzles. Figure~\ref{fig:puzzles} shows Sudoku puzzles, where the complete grid has to be a joint-groups Sudoku. The puzzle on the left is easier. The puzzle on the right contains the minimum possible number of clues: 8. This is because of the following: if there were 7 clues, then at least 2 digits wouldn't have been given, so they could be swapped in the solution.

\begin{figure}[ht!]
\begin{center}
\begin{tikzpicture}
    \TripleEdit
    \FirstColumnTripleEdit{}{}{}{}{}{}{3}{}{5}
    \SecondColumnTripleEdit{}{}{}{}{}{}{}{}{}
    \ThirdColumnTripleEdit{}{}{}{}{}{3}{}{}{}
    \FourthColumnTripleEdit{}{}{}{1}{}{}{}{6}{}
    \FifthColumnTripleEdit{}{}{}{}{}{}{}{}{}
    \SixthColumnTripleEdit{}{}{}{}{}{}{2}{}{}
    \SeventhColumnTripleEdit{7}{}{}{}{}{}{}{}{}
    \EighthColumnTripleEdit{}{}{}{}{}{8}{}{}{6}
    \NinthColumnTripleEdit{}{4}{}{}{5}{}{}{9}{}
\end{tikzpicture} 
\quad \quad
\begin{tikzpicture}
    \TripleEdit
    \FirstColumnTripleEdit{}{6}{}{}{}{}{}{}{}
    \SecondColumnTripleEdit{9}{}{}{}{}{}{}{}{}
    \ThirdColumnTripleEdit{}{}{}{}{}{}{}{}{}
    \FourthColumnTripleEdit{}{}{}{}{}{}{}{}{}
    \FifthColumnTripleEdit{}{}{}{}{}{}{}{}{}
    \SixthColumnTripleEdit{}{}{1}{}{}{}{}{}{}
    \SeventhColumnTripleEdit{5}{7}{8}{}{}{}{}{}{}
    \EighthColumnTripleEdit{}{}{}{}{}{}{}{}{}
    \NinthColumnTripleEdit{}{}{2}{}{}{}{}{}{3}
\end{tikzpicture} 
\end{center}
    \caption{Two JG Sudoku puzzles.}
    \label{fig:puzzles}
\end{figure}

The answers to these puzzles are in Section~\ref{sec:answers}.

\section{An examination of $n=2$}\label{sec:n=2}

We use the knowledge we have gained in the previous sections to completely describe the case of $n=2$.

\subsection{Profiles up to symmetries}

We start by looking at the profiles up to symmetries. We have the following four types:

\begin{enumerate}[A.]
\item Both men prefer the first woman, while both women prefer the first man. 
\item Both men prefer the first woman, while the first woman prefers the first man and the second woman prefers the second man. 
\item The first man prefers the first woman, while the second man prefers the second woman; the first woman prefers the first man, and the second woman prefers the second man.
\item The first man prefers the first woman, while the second man prefers the second woman; the first woman prefers the second man, and the second woman prefers the first man.
\end{enumerate}

\subsubsection{Type A}

The egalitarian costs of the four possible pairs are 2, 3, 3, and 4. Thus we have one pair of soulmates and one hell-pair. Moreover, the soulmates and the hell-pair do not overlap. The presence of soulmates means that there is exactly one stable matching, with one couple being soulmates and the other couple being the hell-couple. The total egalitarian cost is 6.

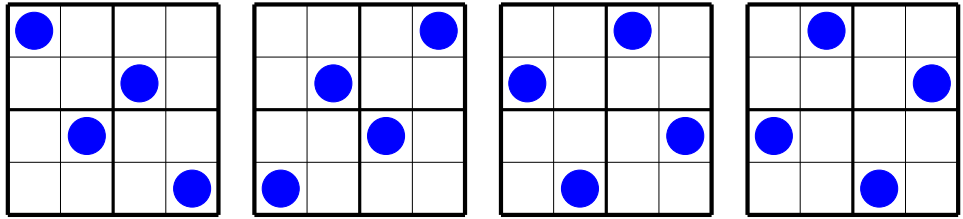
\begin{figure}[ht!]
\begin{center}
\begin{tabular}{cccc}
\centering
    \centering
        \begin{tikzpicture}[scale = 0.7]
            \DoubleEdit
            \filldraw[blue] (-1.5, 1.5) circle (10 pt);
            \filldraw[blue] (-0.5, -0.5) circle (10 pt);
            \filldraw[blue] (0.5, 0.5) circle (10 pt);
            \filldraw[blue] (1.5, -1.5) circle (10 pt);
        \end{tikzpicture} &
        \begin{tikzpicture}[scale = 0.7]
            \DoubleEdit
            \filldraw[blue] (-1.5, -1.5) circle (10 pt);
            \filldraw[blue] (-0.5, 0.5) circle (10 pt);
            \filldraw[blue] (0.5, -0.5) circle (10 pt);
            \filldraw[blue] (1.5, 1.5) circle (10 pt);
        \end{tikzpicture} &
        \begin{tikzpicture}[scale = 0.7]
            \DoubleEdit
            \filldraw[blue] (-1.5, 0.5) circle (10 pt);
            \filldraw[blue] (-0.5, -1.5) circle (10 pt);
            \filldraw[blue] (0.5, 1.5) circle (10 pt);
            \filldraw[blue] (1.5, -0.5) circle (10 pt);
        \end{tikzpicture} &
        \begin{tikzpicture}[scale = 0.7]
            \DoubleEdit
            \filldraw[blue] (-0.5, 1.5) circle (10 pt);
            \filldraw[blue] (1.5, 0.5) circle (10 pt);
            \filldraw[blue] (-1.5, -0.5) circle (10 pt);
            \filldraw[blue] (0.5, -1.5) circle (10 pt);
        \end{tikzpicture} \\
\end{tabular}
\end{center}
    \caption{Profiles of type A.}
    \label{fig:typeA}
\end{figure}

We have four different profiles of type A as shown in Figure~\ref{fig:typeA}. Geometrically, we can describe these profiles as follows. One box has soulmates, and the diagonally opposite box has a hell-pair.

\subsubsection{Type B}

The egalitarian costs of the four possible pairs are 2, 3, 3, and 4. Thus we have one pair of soulmates and one hell-pair. In addition, the soulmates and the hell-pair share a person. The presence of soulmates means that there is exactly one stable matching, with one couple being soulmates and the other couple having an egalitarian cost of 3. The total egalitarian cost is 5.

\begin{figure}[ht!]
\begin{center}
\begin{tabular}{cccc}
\centering
    \centering
        \begin{tikzpicture}[scale = 0.7]
            \DoubleEdit
            \filldraw[blue] (-1.5, 1.5) circle (10 pt);
            \filldraw[blue] (-0.5, -1.5) circle (10 pt);
            \filldraw[blue] (0.5, 0.5) circle (10 pt);
            \filldraw[blue] (1.5, -0.5) circle (10 pt);
        \end{tikzpicture} &
        \begin{tikzpicture}[scale = 0.7]
            \DoubleEdit
            \filldraw[blue] (-1.5, -0.5) circle (10 pt);
            \filldraw[blue] (-0.5, 0.5) circle (10 pt);
            \filldraw[blue] (0.5, -1.5) circle (10 pt);
            \filldraw[blue] (1.5, 1.5) circle (10 pt);
        \end{tikzpicture} &
        \begin{tikzpicture}[scale = 0.7]
            \DoubleEdit
            \filldraw[blue] (-1.5, 0.5) circle (10 pt);
            \filldraw[blue] (-0.5, -0.5) circle (10 pt);
            \filldraw[blue] (0.5, 1.5) circle (10 pt);
            \filldraw[blue] (1.5, -1.5) circle (10 pt);
        \end{tikzpicture} &
        \begin{tikzpicture}[scale = 0.7]
            \DoubleEdit
            \filldraw[blue] (-1.5, -1.5) circle (10 pt);
            \filldraw[blue] (-0.5, 1.5) circle (10 pt);
            \filldraw[blue] (0.5, -0.5) circle (10 pt);
            \filldraw[blue] (1.5, 0.5) circle (10 pt);
        \end{tikzpicture} \\[0.5cm]
        \begin{tikzpicture}[scale = 0.7]
            \DoubleEdit
            \filldraw[blue] (-1.5, 1.5) circle (10 pt);
            \filldraw[blue] (-0.5, -0.5) circle (10 pt);
            \filldraw[blue] (0.5, -1.5) circle (10 pt);
            \filldraw[blue] (1.5, 0.5) circle (10 pt);
        \end{tikzpicture} &
        \begin{tikzpicture}[scale = 0.7]
            \DoubleEdit
            \filldraw[blue] (-1.5, -1.5) circle (10 pt);
            \filldraw[blue] (-0.5, 0.5) circle (10 pt);
            \filldraw[blue] (0.5, 1.5) circle (10 pt);
            \filldraw[blue] (1.5, -0.5) circle (10 pt);
        \end{tikzpicture} &
        \begin{tikzpicture}[scale = 0.7]
            \DoubleEdit
            \filldraw[blue] (-1.5, -0.5) circle (10 pt);
            \filldraw[blue] (-0.5, 1.5) circle (10 pt);
            \filldraw[blue] (0.5, 0.5) circle (10 pt);
            \filldraw[blue] (1.5, -1.5) circle (10 pt);
        \end{tikzpicture} &
        \begin{tikzpicture}[scale = 0.7]
            \DoubleEdit
            \filldraw[blue] (-1.5, 0.5) circle (10 pt);
            \filldraw[blue] (-0.5, -1.5) circle (10 pt);
            \filldraw[blue] (0.5, -0.5) circle (10 pt);
            \filldraw[blue] (1.5, 1.5) circle (10 pt);
        \end{tikzpicture} \\
\end{tabular}
\end{center}
    \caption{Profiles of type B.}
    \label{fig:typeB}
\end{figure}
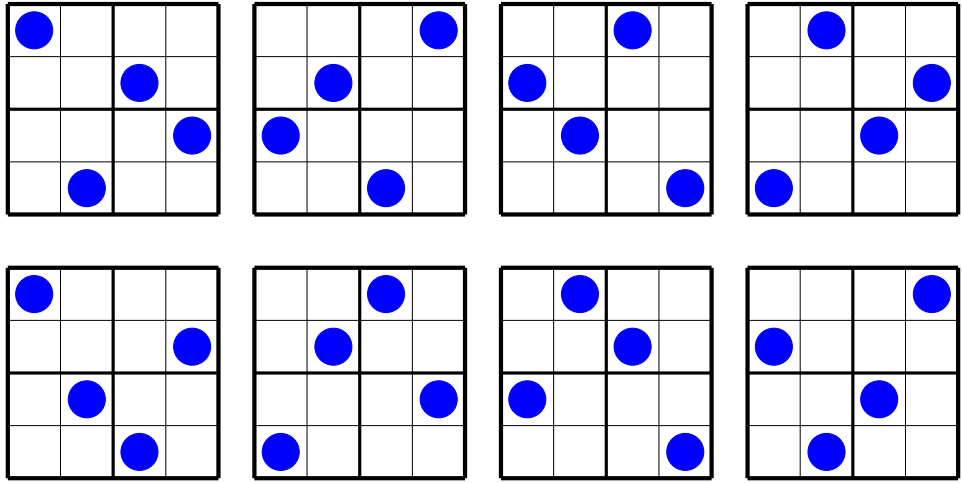

We have eight different profiles of type B as shown in Figure~\ref{fig:typeB}. Geometrically, we can describe them as follows: one box has soulmates, and another box in the same row or column has a hell-pair.

\subsubsection{Type C}

The egalitarian costs of the four possible pairs are 2, 2, 4, and 4. As soulmate pairs can't overlap, the people can be divided into two soulmate pairs or two hell-pairs. The presence of soulmates means that there is exactly one stable matching, in which both couples are soulmates, and the total egalitarian cost is 4.

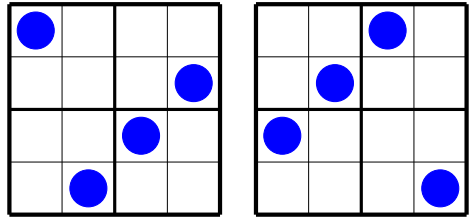
\begin{figure}[ht!]
\begin{center}
\begin{tabular}{cc}
\centering
    \centering
        \begin{tikzpicture}[scale = 0.7]
            \DoubleEdit
            \filldraw[blue] (-1.5, 1.5) circle (10 pt);
            \filldraw[blue] (-0.5, -1.5) circle (10 pt);
            \filldraw[blue] (0.5, -0.5) circle (10 pt);
            \filldraw[blue] (1.5, 0.5) circle (10 pt);
        \end{tikzpicture} &
        \begin{tikzpicture}[scale = 0.7]
            \DoubleEdit
            \filldraw[blue] (-1.5, -0.5) circle (10 pt);
            \filldraw[blue] (-0.5, 0.5) circle (10 pt);
            \filldraw[blue] (0.5, 1.5) circle (10 pt);
            \filldraw[blue] (1.5, -1.5) circle (10 pt);
        \end{tikzpicture}\\
\end{tabular}
\end{center}
    \caption{Profiles of type C.}
    \label{fig:typeC}
\end{figure}

We have two possibilities of type C as shown in Figure~\ref{fig:typeC}. Geometrically, the two pairs of soulmates are in opposite boxes, and the two hell-pairs are in opposite boxes.

\subsubsection{Type D}

The egalitarian costs of the four possible pairs are 3, 3, 3, and 3. There are no soulmates and no hell-pairs. As the only possible mutual egalitarian cost for pairs is $n+1$, this type is a pseudo-Latin profile. There are two possible stable matchings with a total egalitarian cost of 6.

\begin{figure}[h!]
\begin{center}
\begin{tabular}{cc}
\centering
    \centering
        \begin{tikzpicture}[scale = 0.7]
            \DoubleEdit
            \filldraw[blue] (-1.5, 0.5) circle (10 pt);
            \filldraw[blue] (-0.5, -0.5) circle (10 pt);
            \filldraw[blue] (0.5, -1.5) circle (10 pt);
            \filldraw[blue] (1.5, 1.5) circle (10 pt);
        \end{tikzpicture} &
        \begin{tikzpicture}[scale = 0.7]
            \DoubleEdit
            \filldraw[blue] (-1.5, -1.5) circle (10 pt);
            \filldraw[blue] (-0.5, 1.5) circle (10 pt);
            \filldraw[blue] (0.5, 0.5) circle (10 pt);
            \filldraw[blue] (1.5, -0.5) circle (10 pt);
        \end{tikzpicture}\\
\end{tabular}
\end{center}
    \caption{Profiles of type D.}
    \label{fig:typeD}
\end{figure}
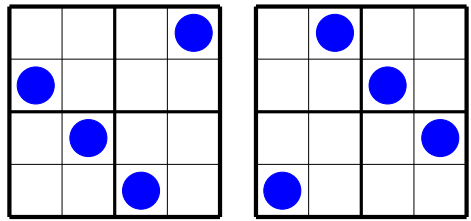

We have two possibilities of type D as shown in Figure~\ref{fig:typeD}. Geometrically, there are no entries in the top-left corner or bottom right corner of each box.

\subsubsection{Summary}

We summarize the results in Table~\ref{tab:sum2}. Types A and B are disjoint profiles because each of the possible mutual rankings is used once, so both of their ranking tally matrices are
$\begin{pmatrix}
1 & 1\\
1 & 1
\end{pmatrix}$.
Types C and D are joint profiles because C has the mutual rankings $(1,1)$ and $(2,2)$, and D has the mutual rankings $(1,2)$ and $(2,1)$. The ranking tally matrices of C and D are
$\begin{pmatrix}2 & 0\\
0 & 2
\end{pmatrix}$
and
$\begin{pmatrix}
0 & 2\\
2 & 0
\end{pmatrix}$
respectively.
In addition, type D is a pseudo-Latin profile because the egalitarian cost of each pair is 3.

\begin{table}[h!]
\begin{center}
\begin{tabular}{|c|c|c|c|c|}
\hline
Type & \# of stable matchings & Total egalitarian cost & \# of profiles & Features\\
\hline
 A & 1 & 6 & 4 & disjoint\\
 B & 1 & 5 & 8 & disjoint\\
 C & 1 & 4 & 2 & joint, mirror\\
 D & 2 & 6 & 2 & joint, pseudo-Latin\\
\hline
\end{tabular}
\end{center}
    \caption{Summary.}
    \label{tab:sum2}
\end{table}

\subsection{Symmetries}

It is interesting to compare profile symmetries with geometric symmetries. 

If we reflect a picture across the main diagonal, we get a picture of the same type. This is because such an operation corresponds to swapping men and women.

The reflection with respect to the middle vertical line corresponds to changing the numbering of men. Namely, man $x$ is swapped with man $n+1-x$. The women's preferences are reversed too: the man ranked $x$ gets new ranking $n+1-x$. We can describe the reflection with respect to the middle horizontal line in a similar manner.

For types A and B, if we reflect a profile across the middle horizontal line or the middle vertical line, we get a profile of the same type, but types C and D swap with each other. 

We see that profile symmetries are very different from geometric symmetries.

\subsection{Combining four non-overlapping profiles together}

Consider a complete Sudoku grid for $n=2$. It corresponds to 4 different preference profiles. As we do not care about the order of the profiles, we consider the Sudoku up to relabeling the digits. 

Equivalently, we can assume that the first row is 1234. The total number of complete Sudokus for $n=2$ is 288. That means the number of grids up to relabeling digits is 12. We present them in lexicographic order. Figure~\ref{fig:3412} shows 4 Sudokus where the second line is 3412 and 2 Sudokus where the second line is 3421.

\begin{figure}[ht!]
\begin{center}
\begin{tabular}{cccc}
\centering
    \centering
        \begin{tikzpicture}[scale = 0.7]
            \DoubleEdit
            \node at (-1.5, 1.5) {\normalsize 1};
            \node at (-1.5, 0.5) {\normalsize 3};
            \node at (-1.5, -0.5) {\normalsize 2};
            \node at (-1.5, -1.5) {\normalsize 4};
            \node at (-0.5, 1.5) {\normalsize 2};
            \node at (-0.5, 0.5) {\normalsize 4};
            \node at (-0.5, -0.5) {\normalsize 1};
            \node at (-0.5, -1.5) {\normalsize 3};
            \node at (0.5, 1.5) {\normalsize 3};
            \node at (0.5, 0.5) {\normalsize 1};
            \node at (0.5, -0.5) {\normalsize 4};
            \node at (0.5, -1.5) {\normalsize 2};
            \node at (1.5, 1.5) {\normalsize 4};
            \node at (1.5, 0.5) {\normalsize 2};
            \node at (1.5, -0.5) {\normalsize 3};
            \node at (1.5, -1.5) {\normalsize 1};
        \end{tikzpicture} &
        \begin{tikzpicture}[scale = 0.7]
            \DoubleEdit
            \node at (-1.5, 1.5) {\normalsize 1};
            \node at (-1.5, 0.5) {\normalsize 3};
            \node at (-1.5, -0.5) {\normalsize 2};
            \node at (-1.5, -1.5) {\normalsize 4};
            \node at (-0.5, 1.5) {\normalsize 2};
            \node at (-0.5, 0.5) {\normalsize 4};
            \node at (-0.5, -0.5) {\normalsize 3};
            \node at (-0.5, -1.5) {\normalsize 1};
            \node at (0.5, 1.5) {\normalsize 3};
            \node at (0.5, 0.5) {\normalsize 1};
            \node at (0.5, -0.5) {\normalsize 4};
            \node at (0.5, -1.5) {\normalsize 2};
            \node at (1.5, 1.5) {\normalsize 4};
            \node at (1.5, 0.5) {\normalsize 2};
            \node at (1.5, -0.5) {\normalsize 1};
            \node at (1.5, -1.5) {\normalsize 3};
        \end{tikzpicture} &
        \begin{tikzpicture}[scale = 0.7]
            \DoubleEdit
            \node at (-1.5, 1.5) {\normalsize 1};
            \node at (-1.5, 0.5) {\normalsize 3};
            \node at (-1.5, -0.5) {\normalsize 4};
            \node at (-1.5, -1.5) {\normalsize 2};
            \node at (-0.5, 1.5) {\normalsize 2};
            \node at (-0.5, 0.5) {\normalsize 4};
            \node at (-0.5, -0.5) {\normalsize 1};
            \node at (-0.5, -1.5) {\normalsize 3};
            \node at (0.5, 1.5) {\normalsize 3};
            \node at (0.5, 0.5) {\normalsize 1};
            \node at (0.5, -0.5) {\normalsize 2};
            \node at (0.5, -1.5) {\normalsize 4};
            \node at (1.5, 1.5) {\normalsize 4};
            \node at (1.5, 0.5) {\normalsize 2};
            \node at (1.5, -0.5) {\normalsize 3};
            \node at (1.5, -1.5) {\normalsize 1};
        \end{tikzpicture} &
        \begin{tikzpicture}[scale = 0.7]
            \DoubleEdit
            \node at (-1.5, 1.5) {\normalsize 1};
            \node at (-1.5, 0.5) {\normalsize 3};
            \node at (-1.5, -0.5) {\normalsize 4};
            \node at (-1.5, -1.5) {\normalsize 2};
            \node at (-0.5, 1.5) {\normalsize 2};
            \node at (-0.5, 0.5) {\normalsize 4};
            \node at (-0.5, -0.5) {\normalsize 3};
            \node at (-0.5, -1.5) {\normalsize 1};
            \node at (0.5, 1.5) {\normalsize 3};
            \node at (0.5, 0.5) {\normalsize 1};
            \node at (0.5, -0.5) {\normalsize 2};
            \node at (0.5, -1.5) {\normalsize 4};
            \node at (1.5, 1.5) {\normalsize 4};
            \node at (1.5, 0.5) {\normalsize 2};
            \node at (1.5, -0.5) {\normalsize 1};
            \node at (1.5, -1.5) {\normalsize 3};
        \end{tikzpicture} 
\end{tabular}
\end{center}
\begin{center}
\begin{tabular}{cc}
\centering
    \centering
        \begin{tikzpicture}[scale = 0.7]
            \DoubleEdit
            \node at (-1.5, 1.5) {\normalsize 1};
            \node at (-1.5, 0.5) {\normalsize 3};
            \node at (-1.5, -0.5) {\normalsize 2};
            \node at (-1.5, -1.5) {\normalsize 4};
            \node at (-0.5, 1.5) {\normalsize 2};
            \node at (-0.5, 0.5) {\normalsize 4};
            \node at (-0.5, -0.5) {\normalsize 1};
            \node at (-0.5, -1.5) {\normalsize 3};
            \node at (0.5, 1.5) {\normalsize 3};
            \node at (0.5, 0.5) {\normalsize 2};
            \node at (0.5, -0.5) {\normalsize 4};
            \node at (0.5, -1.5) {\normalsize 1};
            \node at (1.5, 1.5) {\normalsize 4};
            \node at (1.5, 0.5) {\normalsize 1};
            \node at (1.5, -0.5) {\normalsize 3};
            \node at (1.5, -1.5) {\normalsize 2};
        \end{tikzpicture} &
        \begin{tikzpicture}[scale = 0.7]
            \DoubleEdit
            \node at (-1.5, 1.5) {\normalsize 1};
            \node at (-1.5, 0.5) {\normalsize 3};
            \node at (-1.5, -0.5) {\normalsize 4};
            \node at (-1.5, -1.5) {\normalsize 2};
            \node at (-0.5, 1.5) {\normalsize 2};
            \node at (-0.5, 0.5) {\normalsize 4};
            \node at (-0.5, -0.5) {\normalsize 3};
            \node at (-0.5, -1.5) {\normalsize 1};
            \node at (0.5, 1.5) {\normalsize 3};
            \node at (0.5, 0.5) {\normalsize 2};
            \node at (0.5, -0.5) {\normalsize 1};
            \node at (0.5, -1.5) {\normalsize 4};
            \node at (1.5, 1.5) {\normalsize 4};
            \node at (1.5, 0.5) {\normalsize 1};
            \node at (1.5, -0.5) {\normalsize 2};
            \node at (1.5, -1.5) {\normalsize 3};
        \end{tikzpicture}        
\end{tabular}
\end{center}
    \caption{Sudokus with the second line 3412 and 3421.}
    \label{fig:3412}
\end{figure}

Figure~\ref{fig:4312-4321} shows 2 Sudokus with the second line 4312 and 4 Sudokus with the second line 4321.

\begin{figure}[ht!]
\begin{center}
\begin{tabular}{cc}
\centering
    \centering
        \begin{tikzpicture}[scale = 0.7]
            \DoubleEdit
            \node at (-1.5, 1.5) {\normalsize 1};
            \node at (-1.5, 0.5) {\normalsize 4};
            \node at (-1.5, -0.5) {\normalsize 2};
            \node at (-1.5, -1.5) {\normalsize 3};
            \node at (-0.5, 1.5) {\normalsize 2};
            \node at (-0.5, 0.5) {\normalsize 3};
            \node at (-0.5, -0.5) {\normalsize 1};
            \node at (-0.5, -1.5) {\normalsize 4};
            \node at (0.5, 1.5) {\normalsize 3};
            \node at (0.5, 0.5) {\normalsize 1};
            \node at (0.5, -0.5) {\normalsize 4};
            \node at (0.5, -1.5) {\normalsize 2};
            \node at (1.5, 1.5) {\normalsize 4};
            \node at (1.5, 0.5) {\normalsize 2};
            \node at (1.5, -0.5) {\normalsize 3};
            \node at (1.5, -1.5) {\normalsize 1};
        \end{tikzpicture} &
        \begin{tikzpicture}[scale = 0.7]
            \DoubleEdit
            \node at (-1.5, 1.5) {\normalsize 1};
            \node at (-1.5, 0.5) {\normalsize 4};
            \node at (-1.5, -0.5) {\normalsize 3};
            \node at (-1.5, -1.5) {\normalsize 2};
            \node at (-0.5, 1.5) {\normalsize 2};
            \node at (-0.5, 0.5) {\normalsize 3};
            \node at (-0.5, -0.5) {\normalsize 4};
            \node at (-0.5, -1.5) {\normalsize 1};
            \node at (0.5, 1.5) {\normalsize 3};
            \node at (0.5, 0.5) {\normalsize 1};
            \node at (0.5, -0.5) {\normalsize 2};
            \node at (0.5, -1.5) {\normalsize 4};
            \node at (1.5, 1.5) {\normalsize 4};
            \node at (1.5, 0.5) {\normalsize 2};
            \node at (1.5, -0.5) {\normalsize 1};
            \node at (1.5, -1.5) {\normalsize 3};
        \end{tikzpicture}      
\end{tabular}
\end{center}
\begin{center}
\begin{tabular}{cccc}
\centering
    \centering
        \begin{tikzpicture}[scale = 0.7]
            \DoubleEdit
            \node at (-1.5, 1.5) {\normalsize 1};
            \node at (-1.5, 0.5) {\normalsize 4};
            \node at (-1.5, -0.5) {\normalsize 2};
            \node at (-1.5, -1.5) {\normalsize 3};
            \node at (-0.5, 1.5) {\normalsize 2};
            \node at (-0.5, 0.5) {\normalsize 3};
            \node at (-0.5, -0.5) {\normalsize 1};
            \node at (-0.5, -1.5) {\normalsize 4};
            \node at (0.5, 1.5) {\normalsize 3};
            \node at (0.5, 0.5) {\normalsize 2};
            \node at (0.5, -0.5) {\normalsize 4};
            \node at (0.5, -1.5) {\normalsize 1};
            \node at (1.5, 1.5) {\normalsize 4};
            \node at (1.5, 0.5) {\normalsize 1};
            \node at (1.5, -0.5) {\normalsize 3};
            \node at (1.5, -1.5) {\normalsize 2};
        \end{tikzpicture} &
        \begin{tikzpicture}[scale = 0.7]
            \DoubleEdit
            \node at (-1.5, 1.5) {\normalsize 1};
            \node at (-1.5, 0.5) {\normalsize 4};
            \node at (-1.5, -0.5) {\normalsize 2};
            \node at (-1.5, -1.5) {\normalsize 3};
            \node at (-0.5, 1.5) {\normalsize 2};
            \node at (-0.5, 0.5) {\normalsize 3};
            \node at (-0.5, -0.5) {\normalsize 4};
            \node at (-0.5, -1.5) {\normalsize 1};
            \node at (0.5, 1.5) {\normalsize 3};
            \node at (0.5, 0.5) {\normalsize 2};
            \node at (0.5, -0.5) {\normalsize 1};
            \node at (0.5, -1.5) {\normalsize 4};
            \node at (1.5, 1.5) {\normalsize 4};
            \node at (1.5, 0.5) {\normalsize 1};
            \node at (1.5, -0.5) {\normalsize 3};
            \node at (1.5, -1.5) {\normalsize 2};
        \end{tikzpicture} &
        \begin{tikzpicture}[scale = 0.7]
            \DoubleEdit
            \node at (-1.5, 1.5) {\normalsize 1};
            \node at (-1.5, 0.5) {\normalsize 4};
            \node at (-1.5, -0.5) {\normalsize 3};
            \node at (-1.5, -1.5) {\normalsize 2};
            \node at (-0.5, 1.5) {\normalsize 2};
            \node at (-0.5, 0.5) {\normalsize 3};
            \node at (-0.5, -0.5) {\normalsize 1};
            \node at (-0.5, -1.5) {\normalsize 4};
            \node at (0.5, 1.5) {\normalsize 3};
            \node at (0.5, 0.5) {\normalsize 2};
            \node at (0.5, -0.5) {\normalsize 4};
            \node at (0.5, -1.5) {\normalsize 1};
            \node at (1.5, 1.5) {\normalsize 4};
            \node at (1.5, 0.5) {\normalsize 1};
            \node at (1.5, -0.5) {\normalsize 2};
            \node at (1.5, -1.5) {\normalsize 3};
        \end{tikzpicture} &
        \begin{tikzpicture}[scale = 0.7]
            \DoubleEdit
            \node at (-1.5, 1.5) {\normalsize 1};
            \node at (-1.5, 0.5) {\normalsize 4};
            \node at (-1.5, -0.5) {\normalsize 3};
            \node at (-1.5, -1.5) {\normalsize 2};
            \node at (-0.5, 1.5) {\normalsize 2};
            \node at (-0.5, 0.5) {\normalsize 3};
            \node at (-0.5, -0.5) {\normalsize 4};
            \node at (-0.5, -1.5) {\normalsize 1};
            \node at (0.5, 1.5) {\normalsize 3};
            \node at (0.5, 0.5) {\normalsize 2};
            \node at (0.5, -0.5) {\normalsize 1};
            \node at (0.5, -1.5) {\normalsize 4};
            \node at (1.5, 1.5) {\normalsize 4};
            \node at (1.5, 0.5) {\normalsize 1};
            \node at (1.5, -0.5) {\normalsize 2};
            \node at (1.5, -1.5) {\normalsize 3};
        \end{tikzpicture}       
\end{tabular}
\end{center}
    \caption{Sudokus with the second line 4312 and 4321.}
    \label{fig:4312-4321}
\end{figure}

The Sudokus, in order, have the following sets of profile types:
AAAA, BABA, ABAB, BBBB, BBAA, CDBB, AABB, BBCD, BBBB, CBBD, BDCB, CDCD. There are two BBBB types, but only one for each of the rest. In all 12 of the Sudokus, there are 12 A's, 24 B's, 6 C's, and 6 D's. One can notice that each individual profile appears in exactly 3 of the twelve 2 by 2 Sudokus.

If we order each set of profile types alphabetically, we get the following possibilities. There is one AAAA, four AABBs, two BBBBs, four BBCDs, and one CCDD. We created this list of Sudokus by trying all possibilities. Now, we will explain why these are the only cases.

\begin{enumerate}
\item There are four soulmate pairs and four hell-pairs in each complete grid; so, on average, every profile participating in a Sudoku has one of each. Profiles of types A and B have one of each. Profiles of type C have two of each; profiles of type D have none. It follows that we need the same number of profiles of types C and D in a grid to compensate each other.
\item Every profile of type A overlaps with any profile of type C or D, so A can never be in the same Sudoku as C or D. This explains why any Sudoku with A can only have A's or B's.
\item Each A profile marks either 0 or 2 corners of the Sudoku grid. Each profile of any other type has exactly 1 corner marked. Thus, we need an even number of profiles other than A and consequently an even number of A's.
\end{enumerate}

These three facts taken together show why the only possible Sudokus are AAAA, AABB, BBBB, BBCD, CCDD. The first three types are DG Sudokus, the fifth type is JG Sudoku, and the fourth type is neither. Indeed, the Sudokus can only be disjoint-groups if the profiles are all A's and B's, or joint-groups Sudokus if their profiles are all C's and D's.

The picture of the JG Sudoku is in Figure~\ref{fig:example4profiles}.

\begin{figure}[ht!]
\begin{center}
\begin{tikzpicture}
    \DoubleEdit
    \node at (-1.5, 1.5) {\Large 1};
    \node at (-1.5, 0.5) {\Large 4};
    \node at (-1.5, -0.5) {\Large 3};
    \node at (-1.5, -1.5) {\Large 2};
    \node at (-0.5, 1.5) {\Large 2};
    \node at (-0.5, 0.5) {\Large 3};
    \node at (-0.5, -0.5) {\Large 4};
    \node at (-0.5, -1.5) {\Large 1};
    \node at (0.5, 1.5) {\Large 3};
    \node at (0.5, 0.5) {\Large 2};
    \node at (0.5, -0.5) {\Large 1};
    \node at (0.5, -1.5) {\Large 4};
    \node at (1.5, 1.5) {\Large 4};
    \node at (1.5, 0.5) {\Large 1};
    \node at (1.5, -0.5) {\Large 2};
    \node at (1.5, -1.5) {\Large 3};
\end{tikzpicture}
\end{center}
    \caption{A joint-groups Sudoku example for $n=2$.}
    \label{fig:example4profiles}
\end{figure}

This is the only JG Sudoku for $n=2$ up to relabeling the digits. Notice that this Sudoku has two pairs of identical boxes. The reason is the following: consider the position of digit 1 within the top-left box. There has to be another digit 1 in the exact location relative to its box, but it can't be in the top-right or bottom-left box. So, it has to be in the bottom-right box. This is true for any digit in the top-left box; thus, the bottom-right box has to be a copy of the top-left box. Similarly, the bottom-left and top-right boxes are copies of each other. Thus, all JG Sudokus for $n=2$ are box-cyclic.

The placement matrix for this Sudoku is
\[\begin{pmatrix}
  13 & 24\\ 
  24 & 13
\end{pmatrix}.\]

Profiles 1 and 3 are mirror profiles, and profiles 2 and 4 are pseudo-Latin profiles.

\section{Answers}\label{sec:answers}

Figure~\ref{fig:answers} shows the answers to the JG Sudoku puzzles in Figure~\ref{fig:puzzles} with the given clues highlighted in green.

\begin{figure}[ht!]
\begin{center}
\begin{tikzpicture}
    \filldraw[green] (-9/3,-9/3) rectangle (-7/3,-7/3);
    \filldraw[green] (-9/3,-5/3) rectangle (-7/3,-3/3);
    \filldraw[green] (-5/3,-3/3) rectangle (-3/3,-1/3);
    \filldraw[green] (-3/3,-7/3) rectangle (-1/3,-5/3);
    \filldraw[green] (-3/3,1/3) rectangle (-1/3,3/3);
    \filldraw[green] (1/3,-5/3) rectangle (3/3,-3/3);
    \filldraw[green] (3/3,7/3) rectangle (5/3,9/3);
    \filldraw[green] (5/3,-3/3) rectangle (7/3,-1/3);
    \filldraw[green] (5/3,-9/3) rectangle (7/3,-7/3);
    \filldraw[green] (7/3,5/3) rectangle (9/3,7/3);
    \filldraw[green] (7/3,-1/3) rectangle (9/3,1/3);
    \filldraw[green] (7/3,-7/3) rectangle (9/3,-5/3);
    \TripleEdit
    \FirstColumnTripleEdit{1}{6}{9}{7}{2}{4}{3}{8}{5}
    \SecondColumnTripleEdit{4}{3}{8}{5}{1}{6}{9}{7}{2}
    \ThirdColumnTripleEdit{2}{5}{7}{8}{9}{3}{6}{4}{1}
    \FourthColumnTripleEdit{3}{2}{4}{1}{8}{5}{7}{6}{9}
    \FifthColumnTripleEdit{5}{7}{6}{9}{3}{2}{4}{1}{8}
    \SixthColumnTripleEdit{8}{9}{1}{6}{4}{7}{2}{5}{3}
    \SeventhColumnTripleEdit{7}{8}{5}{3}{6}{9}{1}{2}{4}
    \EighthColumnTripleEdit{9}{1}{2}{4}{7}{8}{5}{3}{6}
    \NinthColumnTripleEdit{6}{4}{3}{2}{5}{1}{8}{9}{7}
\end{tikzpicture} 
\quad \quad
\begin{tikzpicture}
    \filldraw[green] (-7/3,7/3) rectangle (-5/3,3);
   \filldraw[green] (3/3,3/3) rectangle (5/3,3);
    \filldraw[green] (-9/3,5/3) rectangle (-7/3,7/3);
    \filldraw[green] (1/3,3/3) rectangle (3/3,5/3);
    \filldraw[green] (7/3,3/3) rectangle (9/3,5/3);
    \filldraw[green] (7/3,-9/3) rectangle (9/3,-7/3);
    \TripleEdit
    \FirstColumnTripleEdit{1}{6}{4}{5}{7}{8}{2}{3}{9}
    \SecondColumnTripleEdit{9}{2}{5}{3}{1}{6}{8}{4}{7}
    \ThirdColumnTripleEdit{7}{8}{3}{4}{9}{2}{6}{5}{1}
    \FourthColumnTripleEdit{2}{3}{9}{1}{6}{4}{5}{7}{8}
    \FifthColumnTripleEdit{8}{4}{7}{9}{2}{5}{3}{1}{6}
    \SixthColumnTripleEdit{6}{5}{1}{7}{8}{3}{4}{9}{2}
    \SeventhColumnTripleEdit{5}{7}{8}{2}{3}{9}{1}{6}{4}
    \EighthColumnTripleEdit{3}{1}{6}{8}{4}{7}{9}{2}{5}
    \NinthColumnTripleEdit{4}{9}{2}{6}{5}{1}{7}{8}{3}
\end{tikzpicture} 
\end{center}
    \caption{The answers to the JG Sudoku puzzles.}
    \label{fig:answers}
\end{figure}

The corresponding placement matrices are
\[
\begin{pmatrix}
137 & 459 & 268\\ 
268 & 137 & 459\\
459 & 268 & 137
\end{pmatrix}
\quad \textrm{ and } \quad
\begin{pmatrix}
125 & 389 & 467\\ 
367 & 124 & 589\\
489 & 567 & 123
\end{pmatrix}.
\]

Surprisingly, the answer to the first puzzle has a cyclic placement matrix, while the answer to the second puzzle is box-cyclic.

\section{Acknowledgements}

This project was done as part of MIT PRIMES STEP, a program that allows students in grades 6 through 9 to try research in mathematics. Tanya Khovanova is the mentor of this project. We are grateful to PRIMES STEP for this opportunity.

\end{document}